%% file: master.tex
\documentclass[11pt]{article}

\usepackage[margin=1in]{geometry}
\usepackage{amsmath,amsthm,amssymb,amsfonts,float}

\makeatletter
\newcommand{\subjclass}[2][2010]{%
  \let\@oldtitle\@title%
  \gdef\@title{\@oldtitle\footnotetext{#1 \emph{Mathematics subject classification.} #2}}%
}
\makeatother

\usepackage[colorlinks=true, pdfstartview=FitV, linkcolor=blue,citecolor=blue, urlcolor=blue]{hyperref}

\usepackage[abbrev,lite,nobysame]{amsrefs}
\usepackage{times}
\usepackage[usenames,dvipsnames]{color}

\usepackage{esint}

\usepackage{mathtools,enumitem,mathrsfs}

\usepackage[compact]{titlesec}

\usepackage[title]{appendix}

\usepackage{comment}

\mathtoolsset{showonlyrefs=true}
    
\input{macros.tex}

\numberwithin{equation}{section}
    
\begin{document}

\title{Lower bounds on the Lyapunov exponents of stochastic differential equations} 
\author{Jacob Bedrossian\thanks{\footnotesize Department of Mathematics, University of Maryland, College Park, MD 20742, USA \href{mailto:jacob@math.umd.edu}{\texttt{jacob@math.umd.edu}}. J.B. was supported by National Science Foundation CAREER grant DMS-1552826, National Science Foundation RNMS \#1107444 (Ki-Net)} \and Alex Blumenthal\thanks{\footnotesize School of Mathematics, Georgia Institute of Technology, Atlanta, GA 30332, USA \href{mailto:ablumenthal6@gatech.edu}{\texttt{ablumenthal6@gatech.edu}}. A.B. was supported by National Science Foundation grant DMS-2009431} \and Sam Punshon-Smith\thanks{\footnotesize School of Mathematics, Institute for Advanced Study, Princeton, NJ 08540, USA\href{mailto:samuel.punshonsmith@amias.ias.edu}{\texttt{samuel.punshonsmith@amias.ias.edu}}. This material was based upon work supported by the National Science Foundation under Award No. DMS-1803481.}}

\maketitle

\begin{abstract}
  In this article, we review our recently introduced methods for obtaining strictly positive lower bounds on the top Lyapunov exponent of high-dimensional, stochastic differential equations such as the weakly-damped Lorenz-96 (L96) model or Galerkin truncations of the 2d Navier-Stokes equations.
   This hallmark of chaos has long been observed in these models, however, no mathematical proof had been made for either deterministic or stochastic forcing. 

  The method we proposed combines (A) a new identity connecting the Lyapunov exponents to a Fisher information of the stationary measure of the Markov process tracking tangent directions (the so-called ``projective process''); and (B) an $L^1$-based hypoelliptic regularity estimate to show that this (degenerate) Fisher information is an upper bound on some fractional regularity.
  For L96 and GNSE, we then further reduce the lower bound of the top Lyapunov exponent to proving that the projective process satisfies H\"ormander's condition. We review the recent contributions of the first and third author on the verification of this condition for the 2d Galerkin-Navier-Stokes equations in a rectangular, periodic box of any aspect ratio. Finally, we briefly contrast this work with our earlier work on Lagrangian chaos in the stochastic Navier-Stokes equations. We end the review with a discussion of some open problems. 
\end{abstract}

\setcounter{tocdepth}{1}
{\small\tableofcontents}

\section{Lyapunov exponents for stochastic differential equations} \label{sec:Lyap}

%
%



\input{intro.tex}

\section{Formulae for the Lyapunov exponents} \label{sec:formulae}


\input{formulae.tex}

\section{Quantitative lower bounds by the Fisher information}\label{sec:lower-bounds}


\input{qhyp.tex}

\section{Chaos for 2d Galerkin-Navier-Stokes and related models} \label{sec:NSE} 

\input{pspan.tex}

\section{Lagrangian chaos in stochastic Navier-Stokes} \label{sec:LC}

\input{LagChaos.tex}

\section{Looking forward} \label{sec:forward} 

\input{forward.tex}

\addcontentsline{toc}{section}{References}
\bibliographystyle{abbrv}
\bibliography{bibliography}

\end{document}

%% file: macros.tex
\newcommand{\ep}{\epsilon}
\newcommand{\eps}{\epsilon}

\newcommand{\leqc}{\lesssim}

\newcommand{\grad}{\nabla}

\newcommand{\strat}{\circ}
\renewcommand{\div}{\operatorname{div}}

\newcommand{\norm}[1]{\| #1 \|}

\newcommand{\abs}[1]{| #1 |}
\newcommand{\babs}[1]{\left| #1 \right|}
\newcommand{\set}[1]{\{ #1 \}}

\newcommand{\brak}[1]{\left\langle #1 \right\rangle} 

\newcommand{\R}{\mathbb{R}}

\newcommand{\C}{\mathbb{C}}

\newcommand{\Z}{\mathbb{Z}}
\newcommand{\T}{\mathbb{T}}

\renewcommand{\S}{\mathbb{S}}

\newcommand{\loc}{\mathrm{loc}}

\newcommand{\cM}{\mathcal{M}}

\newcommand{\dee}{\mathrm{d}}
\newcommand{\ds}{\dee s}
\newcommand{\dt}{\dee t}

\newcommand{\dx}{\dee x}

\newcommand{\dq}{\dee q}

\DeclareMathOperator{\Div}{\mathrm{div}}
\DeclareMathOperator{\Id}{\mathrm{Id}}
\DeclareMathOperator{\tr}{\mathrm{tr}}

\DeclareMathOperator{\Span}{\mathrm{span}}

\usepackage{bbm}
\newcommand{\1}{\mathbbm{1}}

\renewcommand{\P}{\mathbf{P}}

\newcommand{\E}{\mathbf{E}}
\newcommand{\EE}{\mathbf E}
\newcommand{\PP}{\mathbf P}
\newcommand{\Leb}{\operatorname{Leb}}

\newtheorem{theorem}{Theorem}[section]
\newtheorem{proposition}[theorem]{Proposition}
\newtheorem{corollary}[theorem]{Corollary}

\newtheorem*{lemma*}{Lemma}

\theoremstyle{definition}
\newtheorem{definition}[theorem]{Definition}
\newtheorem{remark}[theorem]{Remark}

%% file: intro.tex

\renewcommand{\phi}{\Phi}


Understanding the ``generic'' long-term dynamics of high (or infinite) dimensional nonlinear systems far from equilibrium remains a daunting task.
In physical applications of interest, many such systems are both subject to unpredictable external forcing and observed to be chaotic in the sense of being very sensitive to the initial condition and forcing. Hence, for all practical purposes, the exact dynamics of any specific trajectory cannot be predicted far in advance and any controlled experiments will not be exactly repeatable. 
 Instead of reckoning such systems one trajectory at a time, a common practice is to view
initial conditions as \emph{random}, i.e., distributed according to some probabilistic law, and 
to attempt to understand how this law evolves as it is transported by the dynamics. 
In this context, the relevant ``time-invariant'' objects are \emph{equilibrium probabilistic laws}
on the phase space of the system, often referred to as \emph{invariant measures} or \emph{stationary measures}.

There is a well-developed abstract theory (smooth ergodic theory) for understanding the invariant measures of chaotic systems, their geometric properties, and how these relate to the asymptotic regimes of trajectories initiated from ``typical'' initial conditions. 
 On the other hand, it is quite hard to verify mathematically that this abstract program applies to systems of practical interest. There are already extremely challenging open problems for vastly simplified 2d toy models of the kinds of chaotic behavior seen in fluid dynamics, e.g., the Chirikov standard map discussed below in Section \ref{subsec:LEchallenges}.

 It turns out that verifying and understanding chaotic properties is far more tractable for systems subjected to \emph{random} noise. 
The kinds of systems we have in mind are, for example, hydrodynamical settings such as with wind over a sail, a weather or climate system, or nonlinear wave systems.
In these settings it has long been suggested to study the random dynamical system generated by the PDE or ODE subjected to random external forcing, and this is often done in applied mathematics (see e.g. \cite{Majda16,BMOV05} and the references therein). 
Even with the simplifications coming from the random forcing, and despite considerable efforts, a thorough, mathematically rigorous, understanding of these random systems is still in its infancy, with many basic open questions remaining.

In this article we will review existing work and our recent contributions \cite{BBPS20,BPS21} in proving that a given system of interest modeled by a stochastic differential equation is chaotic as in high sensitivity to initial conditions for trajectories initiated at Lebesgue-typical points in phase space. The specific systems we apply our methods to are the Lorenz-96 system \cite{Lorenz1996} and Galerkin truncations of the 2d Navier-Stokes equations in a rectangular, periodic box (of any aspect ratio), provided they are subjected to sufficiently strong stochastic forcing\footnote{The deterministic case remains very far out of reach.} (equivalently, sufficiently weak damping) and are sufficiently high dimensional. These are the first results of this type for such models, despite overwhelming numerical evidence (see e.g. \cite{BMOV05,Majda16,KP10,OttEtAl04}). 
Specifically we prove for these models that if the damping parameter is $\eps$, then the top Lyapunov exponent (see Sections \ref{subsec:LEchallenges} and \ref{sec:LEforSDE} for definition) satisfies
\begin{align*}
\lim_{\eps \to 0} \frac{\lambda_1^\eps}{\eps} = \infty
\end{align*}
as $\eps \to 0$ and in particular, $\exists \eps_0 > 0$ such that for all $\eps \in (0, \eps_0)$, $\lambda_1^\eps > 0$. 

\subsubsection*{Outline}

In Section \ref{sec:Lyap} we give background on Lyapunov exponents for stochastic differential equations (SDEs).    
Section \ref{sec:formulae} concerns formulae of Lyapunov exponents through the stationary statistics of tangent directions and contains both classical results and our recent results from \cite{BBPS20} which connects Lyapunov exponents to a certain Fisher information-type quantity.
We discuss in Section \ref{sec:lower-bounds} how to connect the Fisher information to regularity using ideas from hypoellipticity theory (also original work from \cite{BBPS20}), and in Section \ref{sec:NSE} we discuss applications to a class of weakly-driven, weakly-dissipated SDE with bilinear nonlinear drift term (original work in \cite{BBPS20} for Lorenz-96 and for Galerkin Navier-Stokes in \cite{BPS21}). 
In Section \ref{sec:LC} we briefly discuss our earlier related work on Lagrangian chaos in the (infinite-dimensional) stochastic Navier-Stokes equations \cite{BBPS18}.
Finally, in Section \ref{sec:forward} we discuss some open problems and potential directions for research.

\subsection{Lyapunov exponents and their challenges}\label{subsec:LEchallenges}

Let $\phi^t : \R^n \to \R^n, t \in \R_{\geq 0}$ be a flow (autonomous or not) with 
differentiable dependence on initial conditions.
The \emph{Lyapunov exponent} at $x \in \R^n$, when it exists, is the limit
\[
\lambda(x) = \lim_{t \to \infty} \frac1t \log | D_x \phi^t| \, , 
\]
where $D_x \phi^t$ is the Jacobian of $\phi^t$ at $x$, i.e. the derivative with respect to the initial condition. 
Hence, $\lambda(x)$ gives the asymptotic exponential growth rate of the Jacobian as $t \to \infty$. 

The exponent $\lambda(x)$ contains information about the divergence of trajectories: heuristically at least, if $d(x,y)$ is small then 
\[
d(\phi^t(x) , \phi^t(y)) \approx e^{\lambda(x) t} d(x,y)
\]
and hence $\lambda(x) > 0$ implies \emph{exponential sensitivity with respect to initial conditions}, commonly popularized as the ``butterfly effect''. 
Morally, a positive Lyapunov exponent at a `large'  proportion of initial conditions $x \in \R^n$ is a hallmark of chaos, the tendency of a dynamical system to exhibit disordered, unpredictable behavior. In this note we refer to a system such that $\lambda(x) > 0$ for Lebesgue a.e. $x$ as \emph{chaotic}\footnote{We caution the reader that there
 is no single mathematical definition of ``chaos''.  Some definitions refer to the 
 \emph{existence} of a subset of phase space exhibiting chaotic behavior, e.g., Li-Yorke chaos or the presence of a hyperbolic horseshoe. The results discussed in this note pertain to the long-time behavior of Lebesgue-typical initial conditions. }. 

The existence of Lyapunov exponents is usually justified using tools from ergodic theory, and forms a starting point 
for obtaining more refined dynamical features, such as stable/unstable manifolds in the moving frame along `typical' trajectories.
These ideas form the fundamentals of \emph{smooth ergodic theory}, which aims to study \emph{statistical} properties of chaotic systems,
such as decay of correlations, i.e. how $\phi^t(x), t \gg 1$ can `forget' the initial $x \in \R^n$, and 
probabilistic laws such as a strong law of large numbers or central limit theorem for $g \circ \phi^t(x)$, where $g : \R^n \to \R$ is a suitable \emph{observable} of the system; see e.g. discussions in \cite{liu2006smooth,arnold1995random,young2013mathematical,wilkinson2017lyapunov, barreira2002lyapunov}.  

\subsubsection*{A discrete-time example}

Unfortunately, estimating $\lambda(x)$ or proving $\lambda(x) > 0$ for specific systems turns out to be extremely challenging. 
%
A simple, classical model which shows the challenges is the Chirikov standard map family \cite{chirikov1979universal}, written here as 
\[
F_L : \T^2 \circlearrowleft, \quad F_L(x,y) = (2 x + L \sin(2 \pi x) - y, x), 
\]
where $\T^2$ is parametrized as $[0,1)^2$ and both coordinates in $F := F_L$ are taken modulo 1. Here, $L \geq 0$ is a fixed parameter which for purposes of the discussion here will be taken large. The diffeomorphism  $F$ is smooth and volume-preserving, and ergodic theory affirms that the Lyapunov exponent $\lambda(x,y) = \lim_n \frac1n \log |D_{(x,y)} F|$ exists for Leb. a.e. $x$ and satisfies $\lambda(x,y) \geq 0$ where it exists.
The Chirikov standard map itself is frequently used as a toy model of more complicated chaotic systems, e.g., the Navier-Stokes equations in transition from laminar flow to turbulence \cite{mackay1991appraisal}. 

Observe that when $L \gg 1$ and away from an $O(L^{-1})$ neighborhood of $\{ \cos(2 \pi x) = 0\}$, the Jacobian $D_{(x,y)} F$ exhibits strong expansion along tangent directions roughly parallel to the $x$-axis (matched by strong contraction roughly parallel to the $y$-axis). 
In view of this, it is widely conjectured that $\{ \lambda(x) > 0\}$ has positive Lebesgue measure. 
Nevertheless this \emph{standard map conjecture} remains wide open \cite{crovisier2018,pesin2010open}.
A key obstruction is `cone twisting': on long timescales, vectors roughly parallel to the $x$-axis are strongly expanded until the first visit to `critical strip' near $\{ \cos(2 \pi x) = 0\}$, where $DF$ is approximately a rotation by 90 degrees. At this point, vectors roughly parallel to the $x$ axis are rotated to be roughly parallel to the $y$ axis, where strong contraction occurs and previously accumulated expansion can be negated.
Indeed, an estimate on a Lyapunov exponent requires understanding the asymptotic cancellations in the Jacobian as $t \to \infty$.
One manifestation of the subtlety is the wildly tangled coexistence of hyperbolic trajectories \cite{gorodetski2012stochastic} and elliptic islands \cite{duarte1994plenty}.  


The problem of estimating Lyapunov exponents for the standard map is far more tractable in the presence of 
noise / stochastic driving.
Let us consider the standard map subjected to small noise: let $\omega_1, \omega_2, \cdots$ be IID random variables uniformly distributed in $[-\epsilon, \epsilon]$ for some $\epsilon > 0$, and consider the random compositions 
\[
F^n = F_{\omega_n} \circ \cdots \circ F_{\omega_1} \, , \quad F_{\omega_i}(x,y) = F(x + \omega_i, y) \,. 
\]
One can show show that $\forall \epsilon > 0$, the corresponding Lyapunov exponent $\lambda = \lambda(x,y)$ is \emph{deterministic} (independent of the random samples almost surely) and constant (independent of $(x,y)$) with probability 1.
It is a folklore theorem that  $\lambda> 0$ $\forall \epsilon > 0$, while for $L \gg 1$ and $\epsilon \gtrsim e^{-L}$, one can show $\lambda \geq \frac12 \log L$, commensurate with exponential expansion in the $x$-direction over the bulk of phase space \cite{blumenthal2017lyapunov}; in a related vein, see also \cite{shamis2015bounds, blumenthal2020lyapunov, blumenthal2018lyapunov, blumenthal2018positive, lian2012positive}. 

\subsection{Lyapunov exponents for SDE} \label{sec:LEforSDE}
The topic of this note is to discuss developments  in the context of the random dynamical systems generated by stochastic differential equations (SDE), i.e., ODE subjected to Brownian motion driving terms. In this continuous-time framework, numerous additional tools not present in the discrete-time setting become available, e.g., infinitesimal generators, which as we show below, connects the estimation of Lyapunov exponents to regularity estimates (e.g., Sobolev regularity) of solutions to certain (degenerate) elliptic PDE.
A highlight of this approach is our application to the Lyapunov exponents of a class of weakly-driven, weakly-forced SDE,
including famous models such as Lorenz 96 and Galerkin truncations of the Navier-Stokes equations.


For simplicity, in this note we restrict our attention to SDE on $\R^n$, however,  our more general results apply to SDE posed on orientable, geodesically complete, smooth manifolds; see \cite{BBPS20}.  
Let $X_0, X_1, \ldots, X_r : \R^n \to \R^n$ be  smooth vector fields on $\R^n$, and let $W^{1}_t, \ldots, W^{r}_t$ be a collection of independent, real-valued
Brownian motions, with $\Omega$ denoting the corresponding canonical space with probability $\P$ and  $(\mathcal F_t)_{t \geq 0}$ denoting the increasing filtration generated
by $\{W^{k}_s, s \leq t\}_{k = 1}^r$. We consider continuous-time processes $(x_t)$ on $\R^n$
 solving the SDE
\begin{equation}\label{eq:abstractSDE}
\dee x_t = X_0(x_t) \,\dt + \sum_{k =1 }^r X_k(x_t) \circ \dee W_t^{k},
\end{equation}
for fixed initial data $x_0 \in \R^n$. 

Under mild conditions on the vector fields $X_0, \ldots, X_r$ (for example, regularity and the existence of a suitable Lyapunov function to rule out finite time blow-up), global-in-time solutions $(x_t)$ to \eqref{eq:abstractSDE} exist, are unique, and have differentiable dependence of $x_t$ on $x_0$;
in particular, for $\PP$-a.e. $\omega \in \Omega$ and all $t \geq 0$, there exists a stochastic flow of diffeomorphisms $ \phi^t_\omega$ 
 such that $\forall x_0 \in \R^n$, the law of the process $(x_t)_{t \geq 0}$ solving  \eqref{eq:abstractSDE} is the same as that of the process $(\phi_\omega^t(x_0))_{t \geq 0}$; see e.g. \cite{kunita1997stochastic} for the details and general theory of SDEs and stochastic flows.  

 This \emph{stochastic flow of diffeomorphisms} $\phi^t_\omega$ is the analogue of the flow $\phi^t$ corresponding to solutions of the initial value problem of an ODE.
 However, the external stochastic forcing implies a time-inhomogeneity which must be accounted for.
One can show that there exists a $\P$-measure preserving semiflow $\theta^t : \Omega \circlearrowleft, t \geq 0$
corresponding to time-shifts on the Brownian paths, i.e., shifting the path $(W_s)_{s \geq 0}$ to $(W_{t + s} - W_t)_{s \geq 0}$.
Equipped with this time shift, one has the following  with probability 1 and for all $s, t \geq 0$: 
\begin{equation}
\phi^{s + t}_\omega = \phi^t_{\theta^s \omega} \circ \phi^s_\omega.  \label{def:cocycle}
\end{equation}
We now set about summarizing the ergodic theory tools used to study such stochastic flows. 
First, we note that the trajectories $x_t = \phi^t_\omega(x_0)$ for fixed initial 
$x_0 \in \R^n$ form a Markov process adapted to the filtration $(\mathcal F_t)$. 
Moreover,  $\phi^t_\omega$ has \emph{independent increments}: $\forall s, t \geq 0$,  $\phi^s_\omega$ and $\phi^t_{\theta^s \omega}$ are independent.

\subsubsection{Stationary measures and long-term statistics}

\textbf{Markov semigroups. }
We write $P_t(x, A) = \P(\phi^t_\omega(x) \in A)$ for the time-$t$ transition kernel of $(x_t)$. 
Let $\mathcal P_t$ denote the \emph{Markov semigroup} associated to 
 $(x_t)$, defined for bounded, measurable \emph{observables} $h : \R^n \to \R$ by
 \[
 \mathcal P_t h(x) = \E[h(x_t)\,|\,x_0 = x]  = \int_{\R^n} h(y) P_t(x, \dee y). 
 \]
This semigroup gives the expected value of a given observable given a fixed initial condition. 
Via the pairing of functions and measures we derive the (formal) dual $\mathcal{P}_t^\ast$, which gives the evolution of the law of the solution $(x_t)$ given a distribution for the initial condition:
 for a probability measure $\mu_0 \in \mathcal{P}(\R^n)$ and Borel set $A \subseteq \R^n$, 
\begin{align*}
 \mathcal P_t^\ast \mu_0 (A) = \int_{\R^n} P_t(x, A) \dee\mu_0(x) \, .
\end{align*}
That is, $\mathcal P_t^\ast \mu_0$ is the law of $x_t$ assuming $\mu_0$ is the law of $x_0$. 

Taking a time derivative $\partial_t$, we (formally) obtain the \emph{backward Kolmogorov equation} 
\begin{align}
\partial_t \mathcal P_t h(x) = \mathcal L \mathcal P_t h(x) \, , \quad \text{ where } \quad \mathcal L = X_0 + \frac12 \sum_{i = 1}^r X_i^2 \, , \, \label{eq:BKE} 
\end{align}
where, for a given vector field $X$ and $f\in C^\infty$, $Xf$ denotes the derivative of $f$ in the direction $X$. The differential operator $\mathcal{L}$ is called the \emph{(infinitesimal) generator}. 
Assuming that the law of $x_t$ has a density $p_t$ with respect to Lebesgue, the formal dual of \eqref{eq:BKE} is the \emph{Fokker-Plank equation} (or {\em Forward Kolmogorov equation}) given by the following PDE
\begin{align}\label{eq:fokkerPlank}
\partial_t p_t = \mathcal L^* p_t \, , 
\end{align}
where $\mathcal L^*$ denotes the formal $L^2$ adjoint of $\mathcal L$. See e.g. \cite{kunita1997stochastic} for mathematical details.


\medskip

\noindent \textbf{Stationary measures. }
We say a measure $\mu$ is \emph{stationary} if $\mathcal{P}_t^\ast \mu = \mu$. 
That is, if $x_0$ is distributed with law $\mu$, then $x_t$ is distributed with law\footnote{It is important to note that $(x_t)$ itself is \emph{not} constant in $t$; consider e.g. water flowing past a stone in a river.} $\mu$ for all $t > 0$. 
We say that a set $A \subset \R^n$ is \emph{invariant} if $P_t(x, A) = 1$ for all $x \in A$ and $t \geq 0$, 
and we say that a stationary measure $\mu$ is \emph{ergodic} if all invariant sets have $\mu$-measure $0$ or $1$.
By the pointwise ergodic theorem, ergodic stationary measures determine the long-term statistics of a.e. initial datum in their support \cite{DPZ96}:
if $\mu$ is an ergodic stationary measure, then for any bounded, measurable $\varphi:\R^n \to \R$ and $\mu \times \PP$-a.e. $(x,\omega)$ we have that
\[
 \lim_{T \to \infty} \frac{1}{T} \int_0^T \varphi(\Phi^t_\omega(x))\,\dee t = \int_{\R^n} \varphi(x) \dee\mu(x)\,.
\]
Unlike for deterministic systems, stationary measures are usually much easier to characterize for SDEs. In particular, it is often possible to show there exists a unique stationary measure and that it has a smooth density with respect to Lebesgue.
In such a case, Leb. generic initial conditions all have the same long-term statistics, a property often observed in nature and experiments for the physical systems we are interested in.




\medskip

 \noindent \textbf{Existence of stationary measures. }
If the domain of the Markov process were compact (e.g., $\T^n$ instead of $\R^n$) then existence of stationary measures follows from a standard
Krylov-Bogoliubov argument: given an initial probability measure $\mu_0 \in \mathcal{P}(\R^n)$, one considers the time-averaged measures
\[
\bar{\mu}_t := \frac1t \int_0^t \mathcal P_s^* \mu_0 \, \ds.  
\]
Weak-$\ast$ compactness of probability measures on a compact space ensures the sequence $\{\bar{\mu}_t\}_{t \geq 0}$ has a weak-$\ast$ limit point $\mu$ which by construction must be stationary (assuming some mild well-posedness properties for the original SDE). 
On a non-compact domain, one must show tightness of the measures $\set{\bar{\mu}_t}_{t \geq 0}$  (this is essentially saying that solutions do not wander off to infinity too often) and use Prokorov's theorem to pass to the limit in the narrow topology. 
This is often achieved using the method of Lyapunov functions\footnote{These are the probabilistic analogue of Lyapunov's `first method' for ODE, used to ensure convergence to compact attractors. This is not to be confused with Lyapunov exponents, which refer to Lyapunov's `second method'.}/drift conditions \cite{meyn2012markov}, or by using a special structure and the damping in the system (such as the case for e.g. the Navier-Stokes equations \cite{KS}). 

\medskip
\noindent \textbf{Uniqueness of stationary measures. }
The Doob-Khasminskii theorem \cite{DPZ96} implies uniqueness is connected to (A) irreducibility and (B) regularization of the Markov semigroups\footnote{Morally, this is equivalent to how $x \mapsto P_t(x, \cdot)$, i.e., how trajectories with nearby initial conditions have similar statistics.} and in particular, one can deduce that any stationary measure is unique if these properties hold in a sufficiently strong sense. 

Let us first discuss irreducibility.
For a Markov process $(x_t)$ on $\R^n$, we say that $(x_t)$ is \emph{topologically irreducible} if for all open $U \subset \R^n$, $\exists t = t(U, x) \geq 0$ such that
\[
P_t(x, U) > 0 \, .
\]
That is, every initial condition has a positive probability of being in $U$.
This is stronger than necessary to deduce uniqueness, but is sufficient for our discussions.

Regularity is a little more subtle.
A sufficient condition is the requirement of \emph{strong Feller}: 
\begin{align*}
\forall \varphi:\R^n \to \R \quad \textup{ bounded, measurable, } \quad  \mathcal{P}_t \varphi \in C(\R^n;\R), \quad  t > 0.
\end{align*}
For finite dimensional SDEs, it is reasonably common and there exists machinery to characterize this\footnote{In infinite dimensions it is much more rare; luckily it is stronger than what is required just to prove uniqueness (see e.g. \cite{HM06,KNS20}).}.
When $\operatorname{Span}\{ X_i(x) , 1 \leq i \leq r\} = \R^n$ at all $ x\in \R^n$, $\mathcal{L}$ is elliptic and hence strong Feller follows from classical parabolic regularity theory \cite{lieberman1996second} applied to \eqref{eq:BKE} (assuming suitable regularity conditions on the $\set{X_j}$).  
Absent this direct spanning (e.g., when $r < n$), then $\mathcal{L}$ is only degenerate elliptic. However, nearly sharp sufficient conditions for regularization due to $\mathcal{L}$ were derived by H\"ormander \cite{H67}, who obtained a condition (now called \emph{H\"ormander's condition}), in terms of the Lie algebra generated by the vector fields $\{ X_i, 0 \leq i \leq r\}$.
We will return to this important topic of \emph{hypoellipticity} in Section \ref{sec:Hyp}.

\subsubsection{Lyapunov exponents}
We saw that by the ergodic theorem, the long-term behavior of scalar observables is determined by stationary measures. 
A more sophisticated ergodic theorem connects stationary measures to Lyapunov exponents. 
Given $x \in \R^n, v \in \R^n\backslash\{0\}$ ($v$ being considered a direction here) and a random sample $\omega \in \Omega$, the Lyapunov exponent at $(\omega, x, v)$ is defined as the limit (if it exists)
\[
\lambda(\omega, x, v) = \lim_{t \to \infty} \frac{1}{t} \log | D_x \phi^t_\omega v| \,. 
\]
The following (truncated) version of Oseledets' \emph{Multiplicative Ergodic Theorem} (MET) addresses the existence of the limit \cite{raghunathan1979proof,oseledets1968multiplicative,kifer2012ergodic}. 
\begin{theorem}(Oseledets' multiplicative ergodic theorem \cite{oseledets1968multiplicative})\label{thm:MET}
Let $\mu$ be an ergodic stationary measure, and assume a mild integrability condition (see, e.g., \cite{oseledets1968multiplicative,kifer2012ergodic}) then, there exist (deterministic) constants $\lambda_1 > \lambda_2 > \cdots > \lambda_\ell \geq -\infty$ such that for $\P \times \mu$-almost all $(\omega, x) \in \Omega \times \R^n$ and for \emph{all} $v \in \R^n\backslash\{0\}$, the limit defining $\lambda(\omega, x, v)$ exists and takes one of the values $\lambda_i, 1 \leq i \leq \ell$.

Moreover, there exists a $\PP \times \mu$-measurably-varying flag of strictly increasing subspaces 
\[
\emptyset =: F_{\ell+1}(\omega,x) \subset F_{\ell}(\omega,x) \subset \ldots \subset F_1(\omega,x) := \R^n
\]
 such that for $\PP \times \mu$-a.e. $(\omega, x)$ and $\forall v \in F_{j} \setminus F_{j+1}$ there holds
\begin{align*}
\lambda_j = \lim_{t \to \infty} \frac{1}{t} \log | D_x \phi^t_\omega v| = \lambda(\omega,x,v). 
\end{align*}
In particular the top Lyapunov exponent $\lambda_1$ is realized at $\PP \times \mu$-a.e. $(\omega,x)$ and all $v \in \R^n$ outside a positive-codimension subspace $F_2(\omega,x) \subset \R^n$. 
\end{theorem}
We note that under very mild conditions, if the stationary measure $\mu$ is unique, it is 
automatically ergodic; otherwise, each distinct ergodic stationary measure admits its own set
of Lyapunov exponents. 

  The sign of the largest Lyapunov exponent $\lambda_1$ is the most relevant to the stability analysis of typical trajectories, in view of the fact that $\lambda(\omega, x, v) = \lambda_1$ for $v$ in an open and dense set. For this reason we frequently refer to $\lambda_1$ as ``the'' Lyapunov exponent.
The \emph{sum Lyapunov exponent} also turns out to be crucial:
\begin{align*}
\lambda_\Sigma = \sum_{j=1}^\ell m_j \lambda_j = \lim_{t \to \infty} \frac{1}{t} \log | \mathrm{det} \, D_x \phi^t_\omega |, 
\end{align*}
which gives the asymptotic exponential expansion/compression of Lebesgue volume under the flow. Here, $m_j = \dim F_j - \dim F_{j + 1}$ is the \emph{multiplicity} of the $j$-th Lyapunov exponent.

%% file: formulae.tex

Throughout this section, we assume $\phi^t_\omega$ is the stochastic flow of diffeomorphisms corresponding to
the SDE \eqref{eq:abstractSDE} with associated Markov process $x_t = \phi^t_\omega(x), x \in \R^n$. 

\subsection{The projective process}\label{subsec:proj-proc}

As we have seen, Lyapunov exponents are naturally viewed as depending on the tangent direction $v \in \R^n$ at which the derivative $D_x \phi^t_\omega $ is evaluated. For this reason, to estimate Lyapunov exponents it is natural to consider an auxiliary process on \emph{tangent directions} themselves. 
To this end, let $\S \R^n = \R^n \times \S^{n-1}$ denote the unit tangent bundle of $\R^n$, where $\S^{n-1}$ is the unit sphere in $\R^n$. 
Given fixed initial $(x, v) \in \S \R^n$, we define the process $(v_t)$ on $\S^{n-1}$ by
\[
v_t = \frac{D_x \phi^t_\omega(v)}{|D_x \phi^t_\omega(v)|} \,. 
\]
The full process $z_t = (x_t, v_t)$ on $\S \R^n$ is Markovian, and in fact solves an SDE
\[
d z_t = \tilde X_0(z_t) \,\dt + \sum_{i = 1}^r \tilde X_i(z_t) \strat \dee W^{(i)}_t \, , 
\]
where the `lifted' fields $\tilde X_i$ are defined as
\begin{align*}
\tilde X_i(x, v) := \left(X_i(x),  (I - \Pi_v) \nabla X_i(x)v\right). 
\end{align*}
Here, we have written $\Pi_{v} = v \otimes v$ for the orthogonal projection onto the span of $v \in \S^{n-1}$. 
Below, we denote the corresponding generator by 
\[
	\tilde{\mathcal{L}} := \tilde{X}_0 +\frac{1}{2}\sum_{i=1}^r\tilde{X}_i^2.
\] 


\medskip
\noindent {\bf Lyapunov exponents and stationary measures. } Let $(x_t, v_t)$ be a trajectory of the projective process with fixed initial $(x, v) \in \S \R^n$, and observe that at integer times $t \in \mathbb Z_{> 0}$, we have by \eqref{def:cocycle}
\[
\frac{1}{t} \log |D_x \phi^t_\omega(v)| = \frac{1}{t}\sum_{i = 0}^{t - 1} \log |D_{x_i} \phi^1_{\theta^i \omega} v_i|. 
\]
Hence, $\log |D_x \phi^t_\omega|$ is an \emph{additive observable} of $(x_t, v_t)$, i.e., a sum iterated over 
the trajectory $(x_t, v_t)$.
Therefore, the strong law of large numbers for a Markov chain implies the following formula for the Lyapunov exponent: 

\begin{proposition}[See e.g. \cite{kifer2012ergodic}]\label{prop:ergodicThmLE}
Let $\nu$ be an ergodic stationary measure for $(x_t, v_t)$.
Assuming the integral is finite, for $\nu$-a.e. initial $(x, v) \in \S \R^{n-1}$ and $t \geq 0$ we have
\[
t \lambda(\omega, x, v) = \E \int \log |D_x \phi^t_\omega v| \,\dee \nu(x, v) \,. 
\]
with probability 1 ($\E$ denotes integration with respect to $d \P(\omega)$).

Moreover, if $\nu$ is the unique stationary measure for the $(x_t,v_t)$ process, then for $\mu$-a.e. $x$, and \emph{all} $v \in \R^n$, there holds $\lambda_1 = \lambda(\omega,x,v)$ with probability 1 and 
\begin{align}\label{eq:ergodThmFormula}
t \lambda_1 = \E \int \log |D_x \phi^t_\omega v| \,\dee \nu(x, v) \, .
\end{align}
\end{proposition}

\begin{remark}
This latter statement can be interpreted as saying that the existence of a unique stationary measure for the projective process gives a kind of non-degeneracy of the Oseledets' subspace $F_2(\omega,x)$ with respect to $\omega$ \cite{kifer2012ergodic}.
\end{remark}


\medskip
\noindent{\bf A time-infinitesimal version: the Furstenberg-Khasminskii formula. } One of the key benefits
of the SDE framework is the ability to take time derivatives, which turns dynamical questions (e.g., estimates of Lyapunov exponents, identification of stationary densities) into functional-analytic ones (e.g., solutions of degenerate elliptic or parabolic equations) for which many tools are available.
Taking the time derivative of \eqref{eq:ergodThmFormula} gives what is known as the \emph{Furstenberg-Khasminskii formula} (see e.g. \cite{khasminskii2011stochastic,arnold1995random}): 
\begin{proposition} \label{prop:FK}
Assume $(x_t, v_t)$ admits a unique stationary measure $\nu$ on $\S \R^n$ projecting to a stationary measure $\mu$ on $\R^n$ for $(x_t)$.
For $(x, v) \in \S \R^n$, define
\begin{align*}
Q(x) &= \div X_0(x) + \frac12 \sum_{i = 1}^r X_i \div X_i(x) \, ,\\
\tilde Q(x, v) &= \div \tilde X_0(x, v) + \frac12 \sum_{i =1 }^r \tilde X_i \div \tilde X_i(x, v) \,. 
\end{align*}
Then, provided $Q \in L^1(\dee \mu)$ and $\tilde Q \in L^1(\dee \nu)$, there holds 
\begin{gather*}
\lambda_\Sigma = \int Q \,\dee \mu  \quad \text{ and } \\
n \lambda_1 -  \lambda_\Sigma = \int_{\R^n} Q \,\dee \mu - \int_{\S \R^n} \tilde Q \,\dee \nu \,. 
\end{gather*}
\end{proposition}
The first formula expresses $Q(x)$ as the time-infinitesimal rate at which $D_x \phi^t_\omega$ compresses/expands Lebesgue measure, which in this formula is directly related to the asymptotic exponential volume growth/contraction rate $\lambda_\Sigma$. 
Similarly, $\tilde{Q}(x,v)$ is the time-infinitesimal rate at which $D_x \phi^t_\omega$ compresses/expands volume on \emph{the sphere bundle} $\S \R^n = \R^n \times \S^{n-1}$. 
Roughly speaking, contraction of volumes along the $\S^{n-1}$ coordinate is associated  
with expansion in the Jacobian, while expansion of $\S^{n-1}$-volume is related to contraction in the Jacobian; this reversal is the reason for the minus sign in front of $\tilde Q$. 
For some intuition, observe that $(1,0)$ is a sink and $(0,1)$ is a source for the discrete-time system $v_n = A^n v / | A^n v|$ on $S^1$, where $A = \begin{pmatrix} 2 & 0 \\ 0 & 1/2 \end{pmatrix}$. 

\subsection{Sign-definite formulas for Lyapunov exponents}

%


The Furstenberg-Khasminskii formula is highly remarkable in that it reduces the problem of estimating Lyapunov exponents to computing the ensemble average of a single \emph{deterministic} observable, $\tilde Q$, with respect to the stationary measure of $(x_t, v_t)$.
On the other hand, the formula itself is sign-indefinite, as $\tilde Q(x, v)$ takes on both positive and negative values as $(x, v)$ is varied. This is reflective of the cancellation problem mentioned earlier in the estimation of Lyapunov
exponents: tangent growth previously accumulated can be `canceled out' by rotation into 
contracting directions later on in the trajectory.
Hence, without a very precise characterization of $\nu$, it would be very challenging to obtain any useful quantitative estimates on $\lambda_1$ from this formula. 

Given the above, it makes sense to seek a \emph{sign-definite} formula for the Lyapunov exponent.  
Below, given measures $\lambda, \eta, \eta \ll \lambda$ on measurable space $X$, the {\it relative entropy} $H(\eta | \lambda)$ of $\eta$ given $\lambda$ is defined by
\[
H(\eta \,|\, \lambda) = \int_X \log \left( \frac{\dee\eta}{\dee \lambda} \right) \dee \eta \, .
\]
Observe that $H(\eta \,|\, \lambda) \geq 0$, while by strict convexity of $\log$ and Jensen's inequality, we
have $H(\eta | \lambda) = 0$ iff $\eta = \lambda$. 
We also write $\hat \phi^t_\omega : \S \R^d\circlearrowleft$ for the stochastic flow associated to full lifted process $(x_t, v_t)$ on $\S\R^n$; that is, $\hat \phi^t_\omega(x_0, v_0) = (x_t, v_t)$. 
Lastly, given a diffeomorphism $\phi$ of a Riemannian manifold $M$ and a density $g$ on $M$, 
we define $\phi_* g$ to be the density
\[
\phi_* g(x) = g \circ \phi^{-1}(x)  |\det D_x \phi^{-1}| \, ,
\]
noting that if $x$ is distributed like $g \,\dee \mathrm{Vol}_M$, then $\phi(x)$ is distributed like $\phi_* g \,\dee \mathrm{Vol}_M$. 

The following deep formula has its roots in Furstenberg's seminal paper \cite{furstenberg1963noncommuting} and ideas \`a la Furstenberg have been developed by a variety of authors (e.g., \cite{carverhill1987furstenberg,virtser1980products, royer1980croissance, ledrappier1986positivity,baxendale1989lyapunov}), and can be stated as follows: if $\nu \in \mathcal{P}(\S \R^n)$ is a stationary probability measure for the projective process $(x_t, v_t)$ and $\dee\nu(x,v) = \dee\nu_x(v) \dee \mu(x)$ the disintegration of $\nu$, then for all $t > 0$ there holds the following identity (often an inequality in more general settings).

\begin{proposition}[See e.g. \cite{baxendale1989lyapunov} ]\label{prop:relativeEntropyFormulaForLE}
Assume $(x_t, v_t)$ admits a unique stationary measure $\nu$ with density $f = \frac{d \nu}{dq}$, where $\dee q = \dee \mathrm{Vol}_{\S \R^n}$ is Riemannian volume measure on $\S \R^n = \R^n \times \S^{n-1}$. Let $\mu$ be the corresponding stationary measure for $(x_t)$ with density $\rho = \frac{\dee \mu}{\dee x}$. Writing
\[
f_t := (\hat \phi^t_\omega)_* f \,,  \quad \rho_t := (\phi^t_\omega)_* \rho \, , 
\]
we have (under the same integrability condition as Theorem \ref{thm:MET}) 
\[
\E H(\rho_t | \rho) = - t \lambda_\Sigma \, , \quad \text{ and } \quad \E H(f_t | f) = t (n \lambda_1 - 2 \lambda_\Sigma) \, .
\]
\end{proposition}
At least in simple settings, such as for SDEs with a unique stationary measure for the projective process, the formula follows from a slightly more subtle analysis of volume compression/expansion on $\S \R^n$ suitably combined with ergodic theory. Furstenberg \cite{furstenberg1963noncommuting} was the first
to relate relative entropy to Lyapunov exponents; at the generality above, the proof is due to Baxendale \cite{baxendale1989lyapunov}. 

To explore the consequences of Proposition \ref{prop:relativeEntropyFormulaForLE}, let us re-write it into a more suggestive form. 
Let $f_x(v) = f(x,v) / \rho(x), f_{t,x}(v) = f_t(x,v) / \rho_t(x)$ denote the \emph{conditional
densities} of $f$ and $f_t$ along the fiber $\S_x\R^n \simeq \S^{n-1}$. One can then combine the above formulae into the identity
\begin{align}
\E H(f_t|f) - \E H(\rho_t|\rho) = \E \int_{\R^n} H(f_{t,x} | f_x)\dee\mu(x) = t (n \lambda_1 - \lambda_\Sigma)\,. \label{eq:Hhatf}
\end{align}
The left-hand side of this identity is the expectation of a positive quantity, while the right-hand side is non-negative due to the general inequality $n \lambda_1 \geq \lambda_\Sigma$. By the strict convexity, we have
\[
n \lambda_1 = \lambda_\Sigma \quad \iff \quad f_{t,x} \equiv  f_x \text{ with probability 1 for all } t \geq 0 \text{ and }\mu \text{ almost every }x.
\]
Unraveling the definitions, $f_{t,x} \equiv f_x$ means that 
\[
(D_x \phi^t_\omega)_*  f_x =  f_{\phi^t_\omega (x)} \, , 
\]
i.e., the matrices $D_x \phi^t_\omega$, viewed as acting on $\S^{n-1}$ embedded in $\R^n$, transform the conditional density $ f_x$ into the 
density $ f_{\phi^t_\omega(x)}$ of tangent directions at $\phi^t_\omega(x)$. 
This is a very rigid condition in view of the fact that given any two (absolutely continuous) densities $h, h'$ on $\S^{n-1}$, 
\[
\{ A \in GL_n(\R) : A_* h = h'\} 
\]
has empty interior in the space of $n \times n$ matrices.
One can obtain the following beautiful dichotomy by a more detailed analysis of the rigidity in a group of matrices in $SL_n$ that preserve a given probability measure; see, e.g.,  \cite{furstenberg1963noncommuting,baxendale1989lyapunov, ledrappier1986positivity}. 
\begin{theorem}[Furstenberg Criterion] \label{cor:furstCriterion}
Suppose the same setting as Proposition \ref{prop:relativeEntropyFormulaForLE}.  
If $n \lambda_1 = \lambda_\Sigma$, then one of the following holds: 
\begin{itemize}
	\item[(a)] There is a continuously-varying family of inner products $x \mapsto 
	\langle \cdot, \cdot \rangle_x$ with the property that $D_x \phi^t_\omega$ is an isometry
	from $\langle \cdot, \cdot \rangle_{x}$ to $\langle \cdot, \cdot \rangle_{\phi^t_\omega(x)}$ 
	with probability 1 for all $t \geq 0$. 
	\item[(b)] There is a (locally) continuously-varying family of proper subspaces
	$x \mapsto L^i_x \subset \R^d$ with the property that $D_x \phi^t_\omega ( \cup_i L^i_x) = \cup_i L^i_{\phi^t_\omega (x)}$ with probability 1 for all $t \geq 0$. 
\end{itemize}
\end{theorem}
\begin{remark}
Note that in the above, the inner products and the $L^i$ are \emph{deterministic}, which is highly rigid for many random systems. Note that they are also \emph{continuously-varying}.
\end{remark}


However, if one is interested in deducing $\lambda_1> 0$, this criterion is really only useful if $\lambda_\Sigma = 0$, i.e. the system is volume preserving, otherwise one only obtains the non-degeneracy $n\lambda_1 > \lambda_\Sigma$.
Moreover, Theorem \ref{cor:furstCriterion} lacks any quantitative information, and so it cannot be used to obtain concrete estimates with respect to parameters.
Hence, it generally cannot be applied to dissipative systems, even weakly dissipative. 

In the volume preserving case however, criteria \`a la Furstenberg can be a very powerful tool. In our previous work \cite{BBPS18}, we used a suitable (partially) infinite-dimensional extension of  Theorem \ref{cor:furstCriterion} to show that the Lagrangian flow map (i.e. the trajectories of particles in a fluid) is chaotic when the fluid evolves by the stochastically forced 2D Navier-Stokes equations (called \emph{Lagrangian chaos} in the fluid mechanics literature). See Section \ref{sec:LC} for more information. 


\subsection{The best of both worlds: sign-definite and time-infinitesimal}

Proposition \ref{prop:relativeEntropyFormulaForLE} is, on its face, a quantitative and
sign-definite formula for Lyapunov exponents, and this leads to a strong and relatively easy-to-rule-out dichotomy for the degenerate scenario $n \lambda_1 = \lambda_\Sigma$. 
On the other hand, 
the formula itself is not straightforward to work with, requiring both the stationary density 
$f$ for $(x_t, v_t)$ as well as the time-$t$ flow $\phi^t_\omega$ and its derivative $D_x \phi^t_\omega$
as $\omega$ varies. In particular, it is unclear how to glean \emph{quantitative} information beyond the ``soft'' inequality $n\lambda_1 > \lambda_\Sigma$, as would be 
relevant for a damped system (i.e., $\lambda_\Sigma < 0$). 

In view of the sign-indefinite formula \eqref{eq:ergodThmFormula}
and its time-infinitesimal version, the Furstenberg-Khasminskii formula, it is reasonable to hope
that a time-infinitesimal version of Proposition \ref{prop:relativeEntropyFormulaForLE} might exist. 
The authors establish such a formula in our recent work \cite{BBPS20}. 
\begin{proposition}[Theorem A in \cite{BBPS20}]\label{prop:fisherInfo}
Assume $(x_t, v_t)$ has a unique stationary measure $\nu$ with density $f = \frac{\dee \nu}{\dee q}$ on $\S \R^n$. Let $\mu$
denote the corresponding stationary measure for $(x_t)$ on $\R^n$ with density $\rho = \frac{\dee \mu}{\dx}$. 
Define the modified \emph{Fisher informations}
\[
FI(f) = \frac12 \sum_{i =1}^r \int_{\S \R^n} \frac{|\tilde X_i^* f|^2}{f} \,\dee q \, , \quad FI(\rho) = \frac12 \sum_{i =1 }^r \int_{\R^n} \frac{|X_i^* \rho|^2}{\rho} \,\dee x \,. 
\]
Under a mild moment criterion (see \cite{BBPS20}), we have
\[
FI(\rho) = - \lambda_\Sigma \, , \quad \text{ and } \quad FI(f) = n\lambda_1 - 2 \lambda_\Sigma \,. 
\]
Recall that $\tilde X_i^*$ denotes the adjoint of $\tilde X_i$ viewed as an operator on $L^2(dq)$. 
\end{proposition}
\begin{remark}
One can show that $FI(f) - FI(\rho)$ corresponds to an analogous Fisher information on the conditional densities $\hat f_x(v)$, providing the exact time-infinitesimal analogue of \eqref{eq:Hhatf} (see \cite{BBPS20}).  
\end{remark}
These \emph{Fisher-information}-type formulas for Lyapunov exponents enjoy many of the best
qualities of the previous formulas: (A) they are sign-definite, like those in Proposition \ref{prop:relativeEntropyFormulaForLE}, and (B) are also time-infinitesimal like those in Proposition \ref{prop:FK}, and so are inherently simpler, requiring 
only the stationary density $f$ for $(x_t, v_t)$ and how it is acted on by the first-order differential operators $\tilde X_i^*$. 

A key feature of Proposition \ref{prop:fisherInfo} is that a lower bound on $FI(f)$ implies a lower bound on $n\lambda_1 - 2\lambda_\Sigma$. 
The $FI(f)$ itself has the connotation of a \emph{partial regularity} of $f$ along the forcing directions $\tilde X_i$. 
This is reminiscent of techniques in H\"ormander's theory of hypoelliptic operators, where partial
regularity along forcing directions implies regularity in \emph{all} directions under an appropriate Lie algebra spanning condition involving the drift $X_0$. This connection is explored in the next section. 



%% file: qhyp.tex
Let us now set about obtaining quantitative estimates on Lyapunov exponents using the Fisher information as in Proposition \ref{prop:fisherInfo}.
For this, it will be most useful to consider the weakly-forced system 
\begin{align}\label{eq:abstractSDE2}
\dee x_t = X_0^\eps(x_t) \,\dt + \sqrt{\eps}\sum_{k =1 }^r X_k^\eps(x_t) \circ \dee W_t^{k},
\end{align}
where we have also allowed $\eps$ dependence in the vector fields $X_j^\eps$. 
In this case, Proposition \ref{prop:fisherInfo} gives the following Fisher information formula on the stationary density $f^\ep$ of the projective process associated to \eqref{eq:abstractSDE2}
\begin{align}
\frac{\eps}{2} \sum_{j=1}^r \int \frac{\abs{\widetilde{X}^\ast_j f^\eps}^2}{f^\eps} \dee q = n \lambda_1 - 2\lambda_{\Sigma}. 
\end{align}
If $\widetilde{X}_j$ has a bounded divergence\footnote{This is not the case for our examples, but this will not be important as we will eventually work only locally.}, by Cauchy-Schwarz, $\exists C > 0$ such that
\begin{align*}
	\sum_{j=1}^r \norm{\widetilde{X}_j f^\eps}_{L^1}^2 \leq C + FI(f^\eps) = \left(C + \frac{n \lambda_1 - 2\lambda_{\Sigma}}{\eps}\right).
\end{align*}
Hence, we have related  $L^1$-type directional regularity in the forcing directions to the Lyapunov exponents. \emph{If} the lifted forcing directions $\set{\widetilde{X}_j}_{j=1}^r$ spanned the entire tangent space $T_{w} \mathbb S \R^n$ everywhere, then 
we would obtain a lower bound of the Lyapunov exponents of the type 
\begin{align}
\norm{f^\ep}_{\dot{W}^{1,1}}^2 \lesssim \left(1 + \frac{n \lambda_1 - 2\lambda_{\Sigma}}{\eps}\right), \label{ineq:fdotW11}
\end{align}
and so we find a straightforward lower bound on $n\lambda_1 - 2\lambda_{\Sigma}$ in terms of the regularity of $f^\ep$. 
This kind of lower bound is clearly most useful if $\lambda_{\Sigma}$ is small, especially $O(\eps)$, but crucially, it does not have to be \emph{exactly} zero. In this manner, we can treat systems which are close to volume preserving, but not necessarily exactly volume preserving. This is at the crux of why we can treat systems like Lorenz-96 and Galerkin-Navier-Stokes whereas traditional \`a la Furstenberg methods based on e.g. Theorem \ref{cor:furstCriterion}, cannot.    

\subsection{Hypoellipticity} \label{sec:Hyp}

It is not usually the case that $\set{\widetilde{X}_j}_{j=1}^r$ spans $T_{w} \mathbb S \R^n$ and so the lower bound \eqref{ineq:fdotW11} is generally false. 
For example, for additive noise, the lifts satisfy $\widetilde{X}_j = (X_j, 0)$ and so clearly this fails to span $T_w \S\R^n$, regardless of whether or not $\set{X_j}_{j=1}^r$ spans $T_x \R^n$. 
Hence, in general, the Fisher information connects regularity in the lifted forcing directions to the Lyapunov exponents, but a priori, not any other directions in $T_w \mathbb{S} \R^n$. 
For this, we need a concept known as \emph{hypoellipticity}, by which solutions to Kolmogorov equations such as \eqref{eq:BKE} or \eqref{eq:fokkerPlank}, can be smooth even when $\mathcal{L}$ is degenerate, i.e. even when the forcing directions do not span the tangent space.
This effect was studied first by Kolmogorov \cite{K34} in 1934, however clarity on the effect was not fully obtained until H\"ormander's 1967 work \cite{H67}.    

Let us discuss H\"ormander's main insights from \cite{H67}. 
It will make sense to quantify fractional regularity along a vector field $X$ using the group $e^{tX}$ and the $L^p$ H\"older-type semi-norm (brushing aside minor technical details) 
\begin{align*}
\abs{h}_{X,s} := \sup_{t \in (-1,1)} \abs{t}^{-s}\norm{e^{tX}h - h}_{L^p}. 
\end{align*}
H\"ormander's original work was based in $L^2$; our work will be based in $L^1$. For now, we set $p=2$. 

There are two key ideas in \cite{H67}. 
The first, and simpler idea, comes from the Campbell-Baker-Hausdorff formula, which implies for any two vector fields $X,Y$ that (essentially, the Zassenhaus formula):
\begin{align*}
e^{-tX} e^{-tY} e^{tX} e^{tY} = e^{t^2[X,Y] + O(t^3)}, 
\end{align*}
where here $[X,Y]$ is the Lie bracket, i.e. the commutator (see \cite{hochschild1965structure} and \cite{H67}).
In particular, marching forward and then backward by two vector fields $X,Y$ does not quite get back to where it started (unless $X,Y$ commute).  
Therefore we have (using that the $e^{tX}$ are bounded on $L^p$), 
\begin{align*}
\norm{e^{t^2[X,Y] + O(t^3)} - I}_{L^p} \lesssim \norm{e^{tX} - I}_{L^p} + \norm{e^{tY} - I}_{L^p} + \norm{e^{-tX} - I}_{L^p} + \norm{e^{-tY} - I}_{L^p},
\end{align*}
which suggests the remarkable property that any fractional regularity of a function $h$ in directions $X,Y$, i.e. $\abs{h}_{X,s} + \abs{h}_{Y,s} < \infty$, implies $h$ also has (a little less) fractional regularity in the commutator direction $[X,Y]$. 
Another version of Campbell-Baker-Hausdorff (see \cite{H67}) gives
\begin{align*}
e^{t(X+Y)} = e^{tX} e^{tY} e^{t^2[X,Y]}...,
\end{align*}
where the $\ldots$ corresponds to a formal product expansion of higher commutators of $tX$ and $tY$ (and thus higher powers in $t$).
Combined with the previous formal discussion, this suggests that regularity in directions $X,Y$ should also supply regularity in the direction $X+Y$ (and indeed, any linear combination).  
By iterating these heuristics, we get the suggestion that a priori regularity along any set of vector fields $\set{Z_0,...Z_r}$ should imply that there should also be some regularity in \emph{any} direction $Z \in \mathrm{Lie}(Z_0,...,Z_r)$,
where the \emph{Lie algebra} is given by the span of all possible combinations of commutators
\begin{align*}
  \mathrm{Lie}(Z_0,\ldots,Z_r) := \mathrm{span}\set{\mathrm{ad}(Y_m)\ldots\mathrm{ad}(Y_1) Y_0 : Y_j \in \set{Z_0,Z_1,...,Z_r} \, m \geq 0 }, 
\end{align*}
and where $\mathrm{ad}(X)Y := [X,Y]$. In \cite{H67}, these heuristics are made rigorous with the following functional inequality: Suppose that $\forall z \in \R^n$, $\mathrm{Lie}_z(Z_0,...,Z_r) = \set{Z(z): Z \in \mathrm{Lie}(Z_0,...,Z_r)} = T_z \R^n$. Then $\forall s_j \in (0,1)$, $\exists s_\star$ such that for all $0 < s < s_\star$, $\forall R > 0$, and $\forall h \in C^\infty_c(B(0,R))$ there holds 
\begin{align}
\norm{h}_{H^s} \lesssim_{R} \norm{h}_{L^2} +  \sum_{j=0}^r \abs{h}_{Z_j, s_j}. \label{ineq:H1}
\end{align}
In particular, this inequality holds a priori for any  $h \in C^\infty_c(B_R(0))$ and it has nothing to do directly with solutions to any PDE. 
Making this rigorous requires dealing with the errors in the CBH formulas used above.
At any step of the argument, these errors are lower regularity but in new directions, and so dealing with this requires a little finesse and interpolations to close the argument.     

Inequality \eqref{ineq:H1}  is already an interesting observation that can expand the directions of regularity.
In particular, one can use an $L^1$ analogue of \eqref{ineq:H1} to provide a lower bound on the Fisher information based on regularity in any direction contained in the Lie algebra of the forcing directions $\set{\tilde{X}_1,...,\tilde{X}_r}$.
However, H\"ormander was \emph{far} from done. 
Indeed, this is clearly unsatisfying to some degree as this will not even depend on the underlying deterministic dynamical system under consideration, encoded in the drift vector field $\widetilde{X}_0$. Moreover, for additive forcing, \eqref{ineq:H1} fails to add anything at all.
For H\"ormander's second main insight, consider the backward Kolmogorov equation
\begin{align}
\mathcal{L}g = Z_0 g +  \frac{1}{2} \sum_{j=1}^r Z_j^2 g = F. \label{def:ZjK}
\end{align}
Assuming $\set{Z_j}_{j=0}^r$ have bounded divergence\footnote{Alternatively, one can consider the estimates suitably localized.} one obtains the standard $L^2$ ``energy'' estimate:
\begin{align*}
\sum_{j=1}^r\norm{Z_j g}_{L^2}^2 \lesssim \norm{g}_{L^2}^2 + \norm{F}_{L^2}^2. 
\end{align*}
After applying a smooth cutoff $\chi_R(x) = \chi(x/R)$ where $\chi \in C^\infty_c(B_2(0))$, $0 \leq \chi \leq 1$ and $\chi(x) = 1$ for $\abs{x} \leq 1$ and dealing with the commutators as in a Caccioppoli estimate, the functional inequality \eqref{ineq:H1} combined with this estimate implies that \emph{if} $\mathrm{Lie}_z(Z_1,...,Z_r) = T_z \R^n$ at all $z$, then we would obtain an estimate like
\begin{align*}
\norm{\chi_R g}_{H^s} \lesssim_R \norm{g}_{L^2(B_{2R}(0))} + \norm{F}_{L^2(B_{2R}(0))}. 
\end{align*}
However, as discussed above, this condition on the vector fields is often too strong to be useful for us here. 

However, another natural a priori estimate on $g$ is available from \eqref{def:ZjK}.
Indeed, pairing \eqref{def:ZjK} with a test function $\varphi$ we obtain
\begin{align*}
\babs{\int \varphi Z_0 f \dee q} \leq  \frac{1}{2} \sum_{j=1}^r \norm{Z_j^\ast \varphi}_{L^2} \norm{Z_j g}_{L^2} \lesssim \frac{1}{2} \sum_{j=1}^r \left(\norm{\varphi}_{L^2} + \norm{Z_j \varphi}_{L^2} \right) \left(\norm{g}_{L^2} + \norm{F}_{L^2}\right).
\end{align*}
This simple observation shows that for solutions of $\mathcal{L}g = F$, $H^1$-type regularity in the forcing directions automatically provides a corresponding dual $H^{-1}$-type regularity on $Z_0 g$. 
The cornerstone of \cite{H67} is the following functional inequality (i.e. again, not directly related to solutions of any PDEs): 
if $\mathrm{Lie}_z(Z_0,Z_1,...Z_r) = \R^n$ everywhere, then $\exists s \in (0,1)$ such that if $R> 0$ and $h \in C^\infty_c(B_R(0))$, then
\begin{align}
\norm{h}_{H^s} \lesssim \norm{h}_{L^2} + \sup_{\varphi: \norm{\varphi}_{L^2} + \sum_{j=1}^r \norm{Z_j \varphi }_{L^2} \leq 1} \babs{\int \varphi Z_0 h \,\dee q} + \sum_{j=1}^r \norm{Z_j h}_{L^2} =: \norm{h}_{H^1_{\mathrm{hyp}}}. \label{ineq:Hor1} 
\end{align}  
The key heuristic behind this functional inequality is the following observation 
\begin{align}
  \frac{1}{2}\frac{d}{dt} \norm{e^{t Z_0} h - h}_{L^2}^2   & = \brak{e^{t Z_0} h - h, Z_0 e^{tZ_0} h} \\ & \leq \left(\norm{h}_{L^2} + \sum_{j=1}^r \norm{Z_j e^{t Z_0^\ast} (e^{t Z_0}h - h)}_{L^2} \right)\norm{g}_{H^1_{\mathrm{hyp}}}. 
\end{align} 
Therefore, \emph{if} we had something like
\begin{align}
\sum_{j=1}^r \norm{Z_j e^{t Z_0^\ast }(e^{t Z_0}h - h)}_{L^2} \lesssim \sum_{j=1}^r \norm{Z_j h}_{L^2}, \label{ineq:false}
\end{align}
then we could combine the $L^2$ estimate on $\set{Z_j}_{j=1}^r$ with the corresponding dual negative regularity in the $Z_0$ direction to obtain some positive fractional regularity in the $Z_0$ direction, specifically we would have $1/2$ regularity from 
\begin{align*}
\norm{e^{t Z_0} h - h}_{L^2}^2 \lesssim t \norm{h}_{H^1_{\mathrm{hyp}}}^2.
\end{align*}
Unfortunately \eqref{ineq:false} doesn't generally hold\footnote{As in the easier inequalities above, the heuristic \eqref{ineq:false} neglects the creation of higher order commutators, in fact one requires regularity in many other directions in $\mathrm{Lie}(Z_0,...,Z_r)$ as a result.}, and H\"ormander uses a rather ingenious regularization argument to turn this heuristic into reality.  
We shall henceforth call functional inequalities of the type \eqref{ineq:Hor1} \emph{H\"ormander inequalities}. 

The gain in regularity from \eqref{ineq:Hor1} combines with the Kolmogorov equation to get the estimate
\begin{align*}
\norm{g \chi_R}_{H^s} \lesssim \norm{g}_{L^2(B_{2R}(0))} + \norm{F}_{L^2(B_{2R}(0))}, 
\end{align*}
and so provides an analogue of the gain of regularity when studying elliptic equations (though only fractional regularity).
As in that theory, this regularity gain can be iterated to imply that any $L^2$ solution of $\mathcal{L} g=  F$ is $C^\infty$ if $F \in C^\infty$ \cite{H67}.

\subsection{Uniform hypoellipticity}

Next, we want to make the arguments which are quantitative with respect to parameters, and hence we will introduce the notion of \emph{uniform} hypoellipticity.
Let us formalize the definition of H\"ormander's condition for elliptic-type and parabolic-type equations.
For a manifold $M$, we denote $\mathfrak{X}(M)$ the set of smooth vector fields on $M$. 

\begin{definition}[H\"ormander's condition] \label{def:Hormander}
  Given a manifold $\cM$ and a collection of vector fields $\set{Z_0,Z_1,\ldots, Z_r} \subset \mathfrak{X}(\cM)$, we define collections of vector fields $\mathscr{X}_0\subseteq \mathscr{X}_1 \subseteq\ldots$ recursively by
\[
\begin{aligned}
  &\mathscr{X}_0 = \{Z_j\, :\, j \geq 1\},\\
  &\mathscr{X}_{k+1} = \mathscr{X}_k \cup \{[Z_j,Z]\,:\, Z \in \mathscr{X}_k, \quad j \geq 0\}.
  \end{aligned}
\]
We say that $\{Z_i\}_{i=0}^r$ satisfies the {\em parabolic H\"ormander condition} if there exists $k$ such that for all $w \in \mathcal{M}$,
\begin{align}\label{eq:spanningCondIntro22}
    \mathrm{span}\left\{Z(w)\,:\, Z\in \mathscr{X}_k\right\} = T_w \mathcal{M}. 
\end{align}
We say that $\{Z_i\}_{i=0}^r$ satisfies the  {\em (elliptic) H\"ormander condition} if this holds with $\mathscr{X}_0 = \set{Z_j : j \geq 0}$. 
\end{definition}
Note that the parabolic H\"ormander condition is slightly stronger than the elliptic H\"ormander condition. 
\begin{definition}[Uniform H\"ormander's condition] \label{def:UniHormander}
Let $\cM$ be a manifold, and let $\set{Z_0^\epsilon,Z_1^\epsilon,...,Z_r^\epsilon} \subset \mathfrak{X}(\cM)$ be a set of vector fields parameterized by $\epsilon \in (0,1]$.  With $\mathscr{X}_k$ defined as in Definition \ref{def:Hormander} in the parabolic case (resp. elliptic), we say $\set{Z_0^\epsilon,Z_1^\epsilon,...,Z_r^\epsilon}$ satisfies the  uniform parabolic (resp. elliptic) H\"ormander condition on $\cM$ if $\exists k \in \mathbb N$, such that for any open, bounded set $U \subseteq \cM$ there exist constants $\set{K_n}_{n=0}^\infty$, such that for all $\epsilon \in (0,1]$ and all $x \in U$, there is a finite subset $V(x) \subset \mathscr{X}_k$ such that $\forall \xi \in T_x \cM$,
\begin{align*}
\abs{\xi} \leq K_0 \sum_{Z \in V(x)} \abs{Z(x) \cdot \xi} \qquad \sum_{Z \in V(x)}\norm{Z}_{C^n} \leq K_n. 
\end{align*}
\end{definition}
This definition stipulates that any $\eps$ dependence is \emph{locally} (on the manifold) uniform in terms of both regularity and spanning. 
Now we are ready to state the uniform $L^1$-type H\"ormander inequality suitable for use with the Fisher information, proved in \cite{BBPS20}. 
There are many works extending H\"ormander's theory in various ways see e.g., \cite{Polidoro16,Anceschi2019,LanconelliEtAl2020,FarhanGiulio19,GIMV16,Mouhot2018,BL20} and the references therein. However, as far as the authors are aware, there are no works in the $L^1-L^\infty$ framework.
We also  need to consider the  forward Kolmogorov equation $\tilde{\mathcal{L}}^\ast f = 0$, as opposed to the case of the backward Kolmogorov equation considered by H\"ormander \cite{H67}; this changes some details but little of significant consequence is different.

\begin{theorem}[$L^1$-type uniform H\"ormander inequality; Theorem 4.2 \cite{BBPS20}] \label{thm:Hyp}
let $\set{X_0^\eps,X_1^\eps,...,X_r^\eps}$ be a collection of vector fields on $\S \R^n$ satisfying the uniform elliptic H\"ormander condition as in Definition \ref{def:UniHormander}. Then, $\exists s_\star \in (0,1)$ such that if $B_R(x_0) \subset \R^n$ is an open ball and $h \in C^\infty_c(B_R(x_0) \times \S^{n-1})$, then for all $0 < s < s_\star$, $\exists C = C(R,x_0,s)$ such that $\forall \eps \in (0,1)$ there holds the following fractional regularity\footnote{For $s \in (0,1)$ we may define $W^{s,1}$ on a geodesically complete, $n$-dimensional Riemannian manifold with bounded geometry $\cM$ as $$ \norm{w}_{W^{s,1}} = \norm{w}_{L^1} +  \left(\int_{\mathcal{M}} \int_{h \in T_x \mathcal{M} : \abs{h} < \delta_0} \frac{\abs{ w(\mathrm{exp}_x h) - w(x)}}{\abs{h}^{s+n}}\dee h\dee q(x) \right),$$ where $\mathrm{exp}_x: T_x \cM \to \cM$ is the exponential map on $\cM$ and $\dee q$ is the Riemannian volume measure. See e.g. \cite{Triebel} for more details.} estimate \emph{uniformly in $\epsilon$}
\begin{align*}
\norm{h}_{W^{s,1}} \leq C\left(\norm{h}_{L^1} + \sup_{\varphi : \norm{\varphi}_{L^\infty} + \sum_{j=1}^r \norm{X_j^\eps \varphi}_{L^\infty} \leq 1} \babs{ \int \varphi (X_0^\eps)^\ast h\,\dee q} + \sum_{j=1}^r \norm{(X_j^\eps)^\ast h}_{L^1}\right).   
\end{align*}
In particular, applying a smooth cutoff $\chi_R := \chi(x/R)$ for some $\chi \in C^\infty_c(B_2(0))$ with $0 \leq \chi \leq 1$ and $\chi \equiv 1$ if $\abs{x} \leq 1$ to the Kolmogorov equation $\tilde{\mathcal{L}}^\ast f^\eps=0$ (assuming also $\norm{f^\eps}_{L^1} = 1$) and suitably estimating the commutators, we obtain
\begin{equation}\label{eq:lowerbound}
\norm{\chi_R f^\eps}_{W^{s,1}}^2 \lesssim_R 1 + FI(f^\eps). 
\end{equation}
\end{theorem}

\begin{remark}
Hypoellipticity plays a classical role in the theory of SDEs.
In particular, the parabolic H\"ormander condition of Definition \ref{def:Hormander} is exactly the condition most often used to deduce that the Markov semigroup $\mathcal{P}_t$ is strong Feller (the exposition of \cite{H11} is especially intuitive). 
The parabolic H\"ormander condition also often plays a role in proving irreducibility via geometric control theory (see discussions in \cite{Glatt-Holtz2018-wx,Herzog2015-vj,JurdGCT} and specifically in \cite{BBPS20} in regards to the projective process). 
For many applications, it is likely that the parabolic H\"ormander's condition will be used to prove that there exists a unique stationary measure $\nu$ for the projective process (via Doob-Khasminskii \cite{DPZ96}), as required to apply Proposition \ref{prop:fisherInfo}.
Hence the condition of uniformity-in-$\eps$ in Definition \ref{def:UniHormander} will usually be the only additional information required to apply Theorem \ref{thm:Hyp}. 
\end{remark} 

\begin{remark}
Quantitative arguments based on $L^2$ H\"ormander inequalities can be found in \cite{BL20,ABN21} (completed concurrently with or after \cite{BBPS20}).
Thinking about hypoellipticity in terms of functional inequalities, rather than qualitative statements about regularity of solutions to PDEs, has other important advantages as well, for example, it is easier to adapt classical elliptic and parabolic PDE methods, such as De Giorgi or Moser iterations, into hypoelliptic equations \cite{BL20,Mouhot2018,GIMV16}. 
\end{remark}

Obtaining the above Theorem \ref{thm:Hyp} follows an argument generally based on  H\"ormander's original paper \cite{H67}, however, the $L^1-L^\infty$ framework, as opposed to the self-dual $L^2$ framework in \cite{H67}, necessitates a more complicated regularization argument than that used \cite{H67} (which was already quite delicate!).
Moreover, as we are always interested in sphere bundles here, one cannot avoid working on smooth manifolds, which at least under the assumption of geodesic completeness, only adds some technical complexity rather than fundamental difficulties.   

Let us briefly see, heuristically, how one would approach the proof of Theorem \ref{thm:Hyp}.
Motivated by the above discussion regarding \cite{H67}, the main challenge is to obtain $1/2$ of a derivative of $L^1$ H\"older-type regularity in the $\tilde{X}_0^\ast$ ``direction''.
By a bootstrap-type argument, we may assume that we have corresponding regularity along all of the other vector fields in $\mathrm{Lie}_z(\tilde{X}_0,\tilde{X}_1,...,\tilde{X}_r)$ (see \cite{BBPS20} for details).  
Let $S_t$ be a (carefully designed) regularization operator $S_t:L^p \to L^p$. We obtain for any $w \in C^\infty_c(B_R(0) \times \S^{n-1})$, 
\begin{align*}
\norm{e^{t \tilde{X}_0^\ast} w - w}_{L^1} \leq \norm{e^{t \tilde{X}_0^\ast} \left(S_\tau^* w - w\right)}_{L^1} + \norm{S_\tau^* w - w}_{L^1} + \norm{e^{t\tilde{X}_0^*} S_\tau^*w- S_\tau^* w}_{L^1}. 
\end{align*}
We eventually set $\tau \sim \sqrt{t}$ and the regularization operator will be designed so that the first two terms are $O(\tau)$ so we need mainly to work on the latter term, which by duality is estimated by 
\begin{equation*}
	\norm{e^{t\tilde{X}_0^*} S_\tau^*w- S_\tau^* w}_{L^1} \leq  \sup_{\|v\|_{L^\infty} \leq 1} \left| \int_0^t\int_{\S \R^n} (e^{s \tilde{X}_0}v)\, X_0^*S_\tau^*w\, \dq\ds \right|, 
\end{equation*}
and for any fixed $v \in L^\infty$ we have
\begin{align*}
\left|\int_{\S \R^n} (e^{s \tilde{X}_0}v) X_0^* S_\tau^*w\,\dq\right|  & \leq \left|\int_{\S \R^n} (e^{s \tilde{X}_0}v)[\tilde{X}_0,S_\tau]^*w\,\dq\right| + \left|\int_{\S \R^n} (S_\tau e^{s \tilde{X}_0}v)  \tilde{X}_0^* w\,\dq\right|\\
 &\leq \|e^{s X_0}v\|_{L^\infty}\|[\tilde{X}_0,S_\tau]^*w\|_{L^1} + \bigg(\norm{S_\tau e^{s \tilde{X}_0} v}_\infty+ \sum_{j=1}^r\|X_j S_\tau e^{s \tilde{X}_0}v\|_{L^\infty}\bigg) \mathfrak{D}(w), 
\end{align*}
where
\begin{align*}
\mathfrak{D}(h) := \sup_{\varphi : \norm{\varphi}_{L^\infty} + \sum_{j=1}^r \norm{X_j^\eps \varphi}_{L^\infty} \leq 1} \babs{ \int_{\S \R^n} \varphi (\tilde{X}_0^\eps)^\ast h\,\dee q}.
\end{align*}
Hence, the challenge is designing a regularizer such that the commutator $[\tilde{X}_0,S_\tau]^\ast$ loses only $O(\tau^{-1})$ using no a priori regularity in the $\tilde{X}_0$ direction and similarly that $S_\tau$ regularizes the forcing fields $\tilde{X}_j$ like $O(\tau^{-1})$.  
To do this, we let $S_\tau$ be a modified version of H\"ormander's regularizer, which averages the function along directions in $\mathrm{Lie}_z(\tilde{X}_0,...\tilde{X}_r)$ a corresponding amount (higher commutators corresponding to less regularization) in a carefully ordered way.
Specifically, because these `directional mollifiers' do not commute, the order in which they are applied is very important. 
H\"ormander regularized with $S_\tau$, whereas we are fundamentally regularizing with its adjoint $S_\tau^\ast$, which reverses the delicate ordering.
Despite the added difficulty, this turns out to be an important choice for our framework.

%% file: pspan.tex

In this section, we outline how to apply the above ideas to prove a positive Lyapunov exponent for Galerkin truncations of the stochastic 2d Navier-Stokes. 
A general class of models with similar bilinear drift term, which we call \emph{Euler-like} systems, are given by the following SDE: 
\begin{equation}\label{eq:SDEintroFD22}
\dee x_t^\epsilon = (B(x_t^\epsilon,x_t^\eps) - \epsilon A x_t^\epsilon)\dt + \sum_{k=1}^r X_k \dee W_t^{k} \, .
\end{equation}
Here, $\set{X_k}_{k = 1}^r$ is a collection of {\em constant} ($x$-independent) forcing vector fields (i.e. additive forcing) while $B : \R^n \times \R^n \to \R^n$ is a nontrivial (not identically zero) bilinear drift that satisfies
\[
\Div B = 0 \, , \quad x \cdot B(x,x) = 0 \, , 
\]
so in particular the unforced $\epsilon = 0$ dynamics preserve the norm\footnote{In the case of the vorticity form of the 2D Navier-Stokes equations that we will be studying below, this quantity is the \emph{enstrophy}.}, given by $\frac{1}{2}\norm{x}^2$, and volume in $\R^n$ (i.e. the Liouville property).
The term $-\epsilon A$ provides weak linear damping, where $A$ is assumed to be a symmetric, positive-definite $n \times n$ matrix. Stochastically forced versions of the Lorenz 96 model (L96) \cite{Lorenz1996}, Galerkin truncations of 2d and 3d Navier-Stokes on a torus (of arbitrary aspect ratio) \cite{BPS21,Romito2011-ul,E2001-lg} and truncations of commonly used shell models for turbulence \cite{YO87,Gledzer1973,LvovEtAl98,Ditlevsen2010} can be cast in this form. The 2D stochastic Galerkin-Navier-Stokes equations will be described in more detail in Section \ref{sec:SNS-proof} below.

The bilinearity of $B$ implies that solutions can be naturally rescaled into a weakly-damped, weakly-driven system, and the two scalings are equivalent as far as Lyapunov exponents are concerned. Indeed, while the scaling \eqref{eq:SDEintroFD22} is common among models of complex real-world systems, the stationary measure $\mu$ has characteristic energy $\int \abs{x}^2 \dee\mu(x) \approx \eps^{-1}$.
Since we are concerned with the regime $\epsilon \ll 1$, it is natural to rescale and consider a weakly-damped, weakly-driven system. Hence, it is more natural to re-scale so that the long-time behavior remains bounded and non-vanishing as $\ep \to 0$.  By rescaling $x_t^\ep \mapsto \sqrt{\eps} x_{\sqrt{\eps} t}^\eps$, 
replacing $\eps \mapsto \eps^{3/2}$, and using the self-similarity of Brownian motion is equivalent in law to the weakly-driven, weakly damped form
\begin{equation}\label{eq:SDEintroFD22resc}
\dee x_t^\epsilon = (B(x_t^\epsilon,x_t^\eps) - \epsilon A x_t^\epsilon)\dt + \sqrt{\eps}\sum_{k=1}^r X_k \dee W_t^{k}.
\end{equation}
Most importantly, this rescaling does not affect our results on Lyapunov exponents, since upon setting $\hat \eps = \eps^{3/2}$, the Lyapunov exponent $\hat \lambda^{\hat \eps}_1$ of \eqref{eq:SDEintroFD22resc} with parameter $\hat \epsilon$ is related to the Lyapunov exponent $\lambda_1^\eps$ of \eqref{eq:SDEintroFD22} by the identity $\frac{\hat {\lambda}_1^{\hat \eps}}{\hat \eps} = \frac{\lambda_1^\eps}{\eps}$.
This kind of scaling is sometimes called \emph{fluctuation-dissipation} due to the balance between the forcing and the dissipation.


For this class of systems \eqref{eq:SDEintroFD22}, our result below gives a sufficient condition for a positive Lyapunov exponent 
in terms of projective hypoellipticity, i.e., if the lifted vector fields $\{\tilde X_0^\epsilon, \tilde X_1,\ldots \tilde X_r\}$ corresponding to the projective process $(x_t^\epsilon, v_t^\epsilon)$ (denoting $X_0^\epsilon(x) = B(x,x) - \epsilon A x$) satisfy H\"ormander's condition on $\S \R^n$.
\begin{theorem}[Theorem C; \cite{BBPS20}]\label{thm:critForEulerLikeIntro}
Assume that 
\begin{itemize}
\item[(i)] $\{\tilde X_0^\ep,\tilde X_1,\ldots, \tilde X_r\}$ satisfy the elliptic H\"ormander's condition uniformly in $\epsilon \in (0,1)$ as in Definition \ref{def:UniHormander};
\item[(ii)] the  bilinear term $B$ is nontrivial, i.e., $B(x,x) \neq 0$ for some $x \in \R^n$; and
\item[(iii)] the process $(x_t^\eps, v_t^\eps)$ admits a unique stationary density $f^\eps$. 
\end{itemize}
Then, the limit defining the Lyapunov exponent $\lambda_1^\eps$ of \eqref{eq:SDEintroFD22}
exists, and satisfies
$$\lim_{\eps \to 0} \frac{\lambda_1^\eps}{\eps} = \infty \, .$$
In particular, $\exists \eps_0 > 0$ such that for all $\eps \in (0,\eps_0)$ there holds $\lambda_1^\eps > 0$. 
\end{theorem}

A sketch of the proof of Theorem \ref{thm:critForEulerLikeIntro} is given in 
Section \ref{subsec:sketchProofCritEL} below. The most difficult part of applying this result
to a concrete system, e.g., Galerkin-Navier-Stokes, is to prove the parabolic H\"ormander condition for the projective process: general comments on this problem 
are given in Section \ref{subsec:sufficient-cond}, 
while the issue of affirming this for Galerkin-Navier-Stokes is taken up in Section 
\ref{sec:SNS-proof}. 

Given parabolic H\"ormander's condition, unique existence of $f^\eps$ follows, via the 
Doob-Khasminskii theorem, from topological irreducibility of $(x_t^\eps, v_t^\eps)$, i.e., 
the ability to approximately control random trajectories by controlling noise paths. 
For Euler-like models such as \eqref{eq:SDEintroFD22}, this follows from geometric 
control theory arguments and the following well-known cancellation condition on $B(x,x)$
(known to hold for many models such as Galerkin Navier-Stokes, c.f. \cite{Herzog2015-vj,Glatt-Holtz2018-wx}): there exists a collection of vectors $\{e_1,\ldots e_s\}\subset \R^n$ with \[
	\Span\{e_1,\ldots e_s\} = \Span\{X_1,\ldots,X_r\}
	\]
such that for each $1\leq k \leq s$, $B(e_k,e_k) = 0$. For more details see Section 5.3 of \cite{BBPS20}. 
 
\begin{remark}\label{rmk:diffusionTimescale}
The inverse Lyapunov exponent $(\lambda_1^\eps)^{-1}$ is sometimes called the Lyapunov time, 
and is the ``typical'' length of time one must wait for tangent vectors to grow by a factor of $e$. 
Thus, the estimate $\lambda_1^\eps \gg \eps$ implies that the Lyapunov time is $\ll \eps^{-1}$. 
On the other hand, $\eps^{-1}$ is the typical amount of time it takes 
for the Brownian motion $\sqrt{\eps} W_t$ to reach an $O(1)$ magnitude; for this reason it is 
reasonable to refer to $\eps^{-1}$ as a kind of ``diffusion timescale''. So, stated differently, 
our results indicate that as $\eps \to 0$, arbitrarily many Lyapunov times elapse 
before a single ``diffusion time'' has elapsed, indicating a remarkable sensitivity of the Lyapunov exponent to the presence of noise. 

Based on these ideas, one would like to assert that the scaling $\lambda_1^\eps \gg \eps$ 
implies that the deterministic dynamics are ``close'' to positive Lyapunov exponent dynamics, agnostic as to whether the zero-noise system has a positive exponent on a positive area set. However, this assertion does not follow from the scaling $\lambda_1^\eps \gg \eps$ alone: 
even if the Brownian motion itself is small, there could already be a substantial difference between random and corresponding deterministic (zero-noise) trajectories well before time $\eps^{-1}$, e.g., if there is already strong 
 vector growth in the deterministic dynamics. For more on this, see the open problems in Section \ref{sec:forward}. 

\end{remark}

\subsection{Zero-noise limit and rigidity: Proof sketch of Theorem \ref{thm:critForEulerLikeIntro}}
\label{subsec:sketchProofCritEL}
Applying the Fisher information identity (Proposition \ref{prop:fisherInfo}) to the Euler-like system \eqref{eq:SDEintroFD22resc} and using that $\lambda_\Sigma^\ep = -\epsilon \tr A$, we obtain
\begin{align}\label{eq:fisherInfoIdentity11}
FI(f^\ep) = \frac{n \lambda_1^\ep}{\ep} + 2\tr A. 
\end{align}
By the regularity lower bound \eqref{eq:lowerbound}, this implies that for each open ball $B_R(0)$ we have the lower bound
\[
	\|\chi_{R}f^\ep\|_{W^{s,1}}^2 \leqc_R 1 + \frac{\lambda_1^\eps}{\ep},
\]
where the regularity $s \in (0,1)$ and the implicit constant $C = C_R$ are independent of $\epsilon$. 

From this, we see that if $\liminf_\ep \ep^{-1}\lambda_1^\epsilon$ were to remain bounded, then $f^{\epsilon}$ would be bounded in $W^{s,1}_\loc$ uniformly in $\epsilon$. As $W^{s, 1}$ is locally compactly embedded in $L^1$ and $f^\ep$ naturally satisfies certain uniform-in-$\ep$ moment bounds, one can deduce, by sending $\ep\to 0$, that at least one of the following must hold true (see Proposition 6.1, \cite{BBPS20} for details):
\begin{itemize}
    \item[(a)] either $\lim_{\epsilon \to 0} \frac{\lambda_1^\epsilon}{\epsilon} = \infty$; or
    \item[(b)] the zero-noise flow $(x_t^0, v_t^0)$ admits a stationary density $f^0 \in L^1(\S\R^n)$.
\end{itemize}

Let us consider alternative (b). 
While it is natural and common for the projective processes of SDE to admit stationary densities, the existence of an absolutely continuous invariant measure $f^0\dq$ for the projective process of the $\ep=0$ problem
\begin{equation}\label{eq:Euler}
\dot{x}_t = B(x_t,x_t), 
\end{equation}
is quite rigid.  Indeed, in view of the fact that vector growth implies concentration of 
Lebesgue measure in projective space (c.f. the discussion in Section \ref{subsec:proj-proc} after Proposition \ref{prop:FK}), 
the existence of an invariant density 
essentially rules out \emph{any} vector growth for the $\eps = 0$ projective process $(x_t^0, v_t^0)$. 
Precisely, a generalization of Theorem 2.32 in \cite{arnold1999jordan} (see \cite{BBPS20} for details) implies that there is a measurably varying Riemannian metric  $x\mapsto g_x$ such that $\Phi^t$ is an {\em isometry} with respect to $g_x$, namely
\[
	g_x(D_x\Phi^tv, D_x\Phi^tw) = g_{\Phi^t(x)}(v,w),\quad v,w\in T_x\R^n\,,
\]
where $\Phi^t:\R^n\to \R^n$ is the flow associate to the $\ep=0$ dynamics \eqref{eq:Euler}. Hence, we see that if $\liminf_{\eps \to 0}\eps^{-1} \lambda_1 < \infty$, then the deterministic, measure-preserving $\eps=0$ dynamics must be in a situation analogous to possibility (a) in Theorem \ref{cor:furstCriterion}.

In our setting, we show that there is necessarily \emph{some} norm growth as $t \to \infty$ for the $\ep =0$ dynamics due to {\em shearing between conserved energy shells} $\{x\in \R^n\,:\, |x|^2 = E\}$. 
This is straightforward to check: due to the scaling symmetry $\Phi^t(\alpha x) = \alpha\Phi^{\alpha t}(x)$, $\alpha >0$, we have the following orthogonal decomposition of the linearization $D_x\Phi^t$ in the direction $x\in \R^n$:
\[
	D_x\Phi^t x = \Phi^t(x) + t B(\Phi^t(x),\Phi^t(x)) \, , 
\]
noting that $y \cdot B(y,y) \equiv 0$ for all $y \in \R^n$. 
Hence, one obtains the lower bound
\[
	|D_x\Phi^t| \geq t\frac{|B(\Phi^t(x),\Phi^t(x))|}{|x|}
\]
for each $x\in \R^n\backslash\{0\}$ and each $t >0$. 
This contradicts the existence of the Riemannian metric $g_x$ via a Poincar\'e recurrence argument and the fact that the set of stationary points $\{x\in \R^n\,:\, B(x,x) = 0,\, |x|^2 \leq R\}$ is a zero volume set. This is summarized in the following proposition (a proof of which is given in \cite{BBPS20}).

\begin{proposition}[Proposition 6.2 \cite{BBPS20}]\label{prop:mainRigidityIntro}
 Assume that the bilinear mapping $B$ is not identically 0. Let $\nu$ be any invariant probability measure for $\hat \Phi^t$ (the flow corresponding to the (deterministic) $\eps = 0$ projective process) with the property that $\nu(A \times \mathbb S^{n-1}) = \mu(A)$, where $\mu \ll \Leb_{\R^n}$. Then, $\nu$ is singular with respect to volume measure $\dq$ on $\S \R^n$. 
\end{proposition}

\subsection{Verifying projective hypoellipticity: A sufficient condition}\label{subsec:sufficient-cond}

We address here the challenge of verifying the parabolic H\"ormander condition on the sphere bundle $\S\R^n$. Recall that given a smooth vector field $X$ on $\R^n$ we define its lift $\tilde{X}$ to the sphere bundle $\S \R^n$ by
\[
	\tilde{X}(x,v) = \left(X(x),\nabla X(x) v - v\langle v,\nabla X(x)v\rangle\right),
\]
where $\nabla X(x)$ denotes the (covariant) derivative of $X$ at $x$ and is viewed as a linear endomorphism on $T_x\R^n$.
Many of the following general observations about the lifted fields were made in \cite{baxendale1989lyapunov}; see also \cite{BBPS20} for detailed discussions. 

An important property is that the lifting operation can be seen to be a Lie algebra isomorphism onto its range with respect to the Lie bracket, i.e., $[\tilde{X},\tilde{Y}] = [X,Y]{\,\,}^{\widetilde{}}$. Using this observation, the parabolic H\"ormander condition (see Definition \ref{def:Hormander}) on $\S\R^n$ for the lifts of a collection of vector fields $\{X_0,X_1,\ldots,X_r\}\subset \mathfrak{X}(\R^n)$ can be related to non-degeneracy properties of  the Lie sub-algebra $\mathfrak{m}_x(X_0;X_1,\ldots,X_r)$ of $\mathfrak{sl}(T_x\R^n)$ defined by
\begin{equation}\label{eq:m-lie-alg-def}
	\mathfrak{m}_x(X_0;X_1,\ldots,X_r) := \left\{\nabla X(x) - \tfrac{1}{n}\Div X(x) \Id \,:\, X\in \mathrm{Lie}(X_0;X_1,\ldots,X_r)\,,\, X(x) =0\right\},
\end{equation}
where
\[
	\mathrm{Lie}(X_0;X_1,\ldots,X_r) := \mathrm{Lie}(X_1,\ldots,X_r,[X_0,X_1],\ldots,[X_0,X_r]),
\]
is the {\em zero-time ideal} generated by $\{X_0,X_1,\ldots,X_r\}$, with $X_0$ a distinguished ``drift'' vector field (recall that $\mathfrak{sl}_n(T_x \R^n)$ is the Lie algebra of  traceless linear endomorphisms of $T_x\R^n$).

Particularly, if for each $x\in \R^n$, $\mathfrak{m}_x(X_0;X_1,\ldots,X_r)$ acts {\em transitively} on $\S^{n-1}$ in the sense that for each $(x,v)\in \S\R^n$ one has
\begin{equation}\label{eq:transitive}
	\left\{ Av - v\langle v, Av\rangle\,:\, A\in \mathfrak{m}_x(X_0;X_1,\ldots,X_r)\right\} = T_v\S^{n-1},
\end{equation}
then the parabolic H\"ormander condition for $\{X_0,X_1,\ldots,X_r\}$ on $\R^n$ is equivalent to the parabolic H\"ormander condition for the lifts $\{\tilde X_0,\tilde X_1,\ldots,\tilde X_r\}$ on $\S\R^n$. Moreover, the uniform parabolic H\"ormander condition is satisfied on $\S\R^n$ if and only if it is satisfied on $\R^n$ and \eqref{eq:transitive} holds uniformly in the same sense as Definition \ref{def:UniHormander}. 
Since $\mathfrak{sl}(\R^n)$ acts transitively on $\R^n\backslash\{0\}$ (see for instance \cite{boothby1979determination}), a sufficient condition for transitivity on $\S^{n-1}$ is
\[
	\mathfrak{m}_x(X_0;X_1,\ldots X_r) = \mathfrak{sl}(T_x\R^n).
\]

In the specific case of Euler-like models \eqref{eq:SDEintroFD22resc} with $X_0^\ep(x) = B(x,x) -\ep A$ and $\{X_k\}_{k=1}^r$ as in \eqref{eq:SDEintroFD22resc}, the situation can be simplified if $\mathrm{Lie}(X_0;X_1,\ldots,X_r)$ contains the constant vector fields $\{\partial_{x_k}\}_{k=1}^n$. In this case, the family of $x$ and $\eps$-independent endomorphisms
\[
	H_k := \nabla [\partial_{x_k},X_0^\ep] = \nabla [\partial_{x_k},B],\quad k=1,\ldots n,
\]
generate the Lie algebra $\mathfrak{m}_x(X_0^\ep;X_1,\ldots,X_r)$ at all $x \in \R^n$.
This argument implies the following sufficient condition for projective spanning.  
\begin{corollary}[See \cite{BBPS20}]
Consider the bilinear Euler-like models \eqref{eq:SDEintroFD22resc}. If $\mathrm{Lie}(X_0;X_1,\ldots,X_r)$ contains $\{\partial_{x_k}\}_{k=1}^n$, then $\{\tilde{X}_0^\ep,\tilde X_1, \ldots \tilde X_r\}$ satisfy the uniform parabolic H\"ormander condition (in the sense of Definition \ref{def:UniHormander}) on $\S\R^n$ if
\begin{equation}\label{eq:sufficient condition}
	\mathrm{Lie}(H^1,\ldots, H^n) = \mathfrak{sl}(\R^n).
\end{equation}
\end{corollary}
\noindent This criterion is highly useful, having reduced projective spanning to a question about a single Lie algebra of trace-free matrices.

In \cite{BBPS20}, we verified this condition directly for the Lorenz 96 system \cite{Lorenz1996}, which is defined for $n$ unknowns in a periodic array  by the nonlinearity $B$ given by 
\begin{align}
B_\ell(x,x) = x_{\ell+1}x_{\ell -1} - x_{\ell-2}x_{\ell-1}.  \label{def:L96B}
\end{align}
The traditional case is $n=40$, but it can be considered in any finite dimension. 
In particular we proved the following.
\begin{corollary}[Corollary D; \cite{BBPS20}]
Consider the L96 system given by \eqref{eq:SDEintroFD22resc} with the nonlinearity \eqref{def:L96B} and $X_k = q_k e_k$ for $k \in \set{1,..,r}$, $q_k \in \R$, and $e_k$ the canonical unit vectors. If $q_1, q_2 \neq 0$ and $n \geq 7$, then
\begin{align*}
\lim_{\eps \to 0} \frac{\lambda_1^\eps}{\eps} = \infty \, .
\end{align*}
In particular $\exists \eps_0 > 0$ such that $\lambda_1^\eps > 0$ if $\eps \in (0, \eps_0)$. 
\end{corollary}

\subsection{Projective hypoellipticity for 2d Galerkin-Navier-Stokes}\label{sec:SNS-proof}
Let's now see how we can go about verifying the projective hypoellipticity condition for a high-dimensional model of physical importance, namely Galerkin truncations of the 2d stochastic Navier-Stokes equations on the torus of arbitrary side-length ratio $\T^2_r = [0,2\pi) \times [0, \frac{2\pi}{r})$ (periodized) for $r > 0$. Recall that the Navier-Stokes equations on $\T^2_r$ in vorticity form are given by 
\begin{equation}
\partial_t w + u \cdot \grad w - \eps \Delta w = \sqrt{\eps} \dot{W}_t,  
\end{equation}
where $u$ is the divergence free velocity field satisfying the Biot-Savart law $u = \nabla^{\perp}(-\Delta)^{-1}w$ and $\dot{W}_t$ is a white-in time, colored-in-space Gaussian forcing which we will take to be diagonalizable with respect to the Fourier basis with Fourier transform supported on a small number of modes.

In the work \cite{BPS21} by the first and last authors of this note, we consider a Galerkin truncation of the 2d stochastic Navier-Stokes equations at an arbitrary frequency $N\geq 1$ in Fourier space by projecting onto the Fourier modes in the truncated lattice
\[
	\Z_{0,N}^2 := \{(k_1,k_2) \in \Z^{2}\backslash\{0\}\,:\, \max\{|k_1|,|k_2|\}\leq N\}\subseteq \Z^2,
\]
giving rise to a $n=|\Z^2_{0,N}| = (2N+1)^2-1$ dimensional stochastic differential equation with the reality constraint $w_{-k} = \overline{w}_k$ for $w = (w_k)\in \C^{\Z^{2}_{0,N}}$ (that is, the vector is indexed over $\Z^2_{0,N}$) governed by
\begin{equation}
	\dee w_k = (B_k(w,w) - \ep \abs{k}^2_r w_k)\dt + \sqrt{\ep}\dee W^k,\label{def:NSEFF}
\end{equation}
where $\abs{k}_r^2 := k_1^2 + r^2 k_2^2$, and $W^k_t = \alpha_k W^{a,k}_t + i\beta_k W^{b,k}_t$ are independent complex Wiener processes satisfying $W^k_t = \overline{W}^{-k}_t$ ($W^{a,k}_t,W^{b,k}_t$ are standard iid Wiener processes) with $\alpha_k$, $\beta_k$ arbitrary such that $\alpha_k = 0 \Leftrightarrow \beta_k= 0$. The symmetrized non-linearity $B_k(w,w)$ is given by
\begin{equation}\label{eq:B-def}
	B_k(w,w) := \frac{1}{2}\sum_{j+\ell=k} c_{j,\ell} w_j w_\ell,\quad c_{j,\ell} := \langle j^\perp,\ell\rangle_r\left(\frac{1}{|\ell|^2_r} - \frac{1}{|j|^2_r}\right)
\end{equation}
where the sum runs over all $j,\ell \in Z^2_{0,N}$ such that $j+\ell=k$ and we are using the notation $\langle j^\perp,\ell\rangle_r := r (j_2\ell_1 - j_1\ell_2)$. In what follows the coefficient $c_{j,\ell}$ always depends on $r$ but we suppress the dependence for notational simplicity.

We will regard the configuration space $\C^{\Z^{2}_{0,N}}$ as a complex manifold with complexified tangent space spanned by the complex basis vectors $\{\partial_{w_k}\,:\,k\in \Z^2_{0,N}\}$ (Wirtinger derivatives) satisfying $\partial_{w_{-k}} = \bar{\partial}_{w_k}$. See \cite{huybrechts2005complex} for the notion of complexified tangent space and \cite{BPS21} for discussion on how to use this complex framework for checking H\"ormander's condition.
In this basis, we can formulate the SDE \eqref{def:NSEFF} in the canonical form
\begin{equation}\label{eq:SDE-SNS-canon}
	\dee w_t = X_0^\ep(w_t) + \sum_{k\in \mathcal{Z}^0} \sqrt{\ep}\partial_{w_k}\dee W^k_t,
\end{equation}
where the drift vector field $X_0^\ep$ is given by $X_0^\ep(w) := \sum_{k\in \Z^{2}_{0,N}} (B_k(w,w) - \ep|k|^2_rw_k)\partial_{w_k}$ and the set of driving modes $\mathcal{Z}^0$ is given by $\mathcal{Z}^0 := \{k\in \Z^2_{0,N} \,:\, \alpha_k, \beta_k \neq 0\}$.

As in the setting of \cite{E2001-lg,HM06}, we consider very degenerate forcing and study how it spreads throughout the system via the nonlinearity $B_\ell(w,w)$. Specifically, we define the sets 
\[
\mathcal{Z}^n = \{\ell \in \Z^2_{0,N}\,:\,\ell = j+k,\, j\in \mathcal{Z}^0,\, k\in \mathcal{Z}^{n-1}, c_{j,k}\neq 0 \}, \quad n\geq 0
\]
and assume that the driving modes $\mathcal{Z}^0$ satisfy $\bigcup_{n\geq 0} \mathcal{Z}^n = \Z^2_{0,N}$. Under this assumption on $\mathcal{Z}^0$ it can be shown (see \cite{BBPS20} Proposition 3.6 or \cite{E2001-lg,HM06}) that the complexified Lie algebra $\mathrm{Lie}(X^\ep_0;\{\partial_{w_k}\,:\,k\in \mathcal{Z}^0\})$ contains the constant vector fields $\{\partial_{w_k}\,:\,k\in \Z^2_{0,N}\}$ and therefore satisfies the uniform parabolic H\"ormander condition on $\C^{\Z^{2}_{0,N}}$.

\subsubsection{A distinctness condition on a diagonal subalgebra}

As discussed in Section \ref{subsec:sufficient-cond}, in order to verify projective hypoellipticity for the vector fields $X^\ep_0;\{\partial_{w_k}\,:\,k\in \mathcal{Z}^0\}$, it suffices to study the generating properties of a suitable matrix Lie algebra. 
In \cite{BPS21}, we show this can be reformulated to a condition on the constant, {\em real valued} matrices $H^k = \nabla[\partial_{w_k},B]$, $k\in \Z^{2}_{0,N}$, represented in $\{\partial_{w_k}\}$ coordinates by $(H^k)_{\ell,j} = \partial_{w_j}\partial_{w_k}B_\ell(w,w) = c_{j,k}\delta_{\ell=j+k}$. After obtaining this reformulation, the main result of \cite{BPS21} is the following non-degeneracy property of the matrices $\{H^k\}$.
\begin{theorem}[Theorem 2.13, \cite{BPS21} (see also Proposition 3.11)]\label{thm:2d-GSNS-Lie-gen}
Consider the 2d stochastic Galerkin Navier-Stokes equations with frequency truncation $N$ on $\T_r^2$
and suppose that $N\geq 392$. Then, the following holds: 
\begin{equation}\label{eq:Lie-suff-H}
\mathrm{Lie}(\{H^k:k\in \Z^2_{0,N}\}) = \mathfrak{sl}_{\Z^2_{0,N}}(\R) \, ,
\end{equation}
where $\mathfrak{sl}_{\Z^2_{0,N}}(\R)$ denotes the Lie algebra of real-valued traceless matrices indexed by the truncated lattice $\Z^2_{0,N}$. Therefore projective hypoellipticity holds for \eqref{eq:SDEintroFD22resc} and by Theorem \ref{thm:critForEulerLikeIntro} the top Lyapunov exponent satisfies $\lim_{\ep\to 0}\ep^{-1}\lambda_1^{\ep} = \infty$.
 \end{theorem}

\begin{remark}
Verifying the Lie algebra generating condition \eqref{eq:Lie-suff-H} can be quite challenging due to the fact that there there are $n = |\Z^2_{0,N}|$ matrices and $n^2 - 1$ degrees of freedom to span.
The matrices are also banded in the sense that for each $k$, $(H^k)_{\ell,j}$ couples most of the lattice values $\ell, j$ along the band $k=\ell-j$ and therefore it is extremely challenging to isolate elementary matrices (matrices with only one non-zero entry) as one can do in ``local in frequency'' models like L96 \eqref{def:L96B} (see \cite{BBPS20}).
Moreover, brute force computational approaches that successively generate Lie bracket generations and count the rank by Gaussian elimination (such as the Lie-Tree algorithm in \cite{Elliott2009-ha}) are only available for fixed $r \in \R_+$ and $N \in \Z_+$, and can be subject to numerical error (for instance if $r$ is chosen irrational) which destroy the validity of the proof.
\end{remark}

In order to show that \eqref{eq:Lie-suff-H} holds, in \cite{BPS21} we take an approach inspired by the root-space decomposition of semi-simple Lie algebras and study genericity properties of the following diagonal sub-algebra of $\mathrm{Lie}(\{H^k\})$
\[
	\mathfrak{h} := \Span\{\mathbb{D}^k\,:\,k\in \Z^{2}_{0,N}\}, 
\]
where $\mathbb{D}^k = [H^k,H^{-k}]$ are a family of diagonal matrices with diagonal elements $\mathbb{D}_i^k = (\mathbb{D}^k)_{ii}$ given by
\[
	\mathbb{D}_i^k = c_{i,k}c_{i+k,k}\1_{\Z^2_{0,N}}(i+k) - c_{i,k}c_{i-k,k}\1_{\Z^2_{0,N}}(i-k).
\]

Using that, for a given diagonal matrix $\mathbb{D}\in \mathfrak{sl}_{\Z^2_{0,N}}(\R)$, the adjoint action $\mathrm{ad}(\mathbb{D}):\mathfrak{sl}_{\Z^2_{0,N}}(\R)\to \mathfrak{sl}_{\Z^2_{0,N}}(\R)$, where $\mathrm{ad}(\mathbb{D})H = [\mathbb{D},H]$, has eigenvectors given by the elementary matrices $E^{i,j}$ (i.e. a matrix with a
one in the $i$th row and $j$th column and zero elsewhere) $\mathrm{ad}(\mathbb{D})E^{i,j} = (\mathbb{D}_i -\mathbb{D}_j)E^{i,j}$, means that $\mathrm{ad}(\mathbb{D})$ has a simple spectrum if the diagonal entries of $\mathbb{D}$ have {\em distinct differences}, $\mathbb{D}_i - \mathbb{D}_j \neq \mathbb{D}_{i^\prime} - \mathbb{D}_{j^\prime}$, $(i,j)\neq (i^\prime, j^\prime)$. This implies that if $H$ is a matrix with non-zero non-diagonal entries and $\mathbb{D}$ has distinct differences, then for $M= n^2 - n$, the Krylov subspace 
\[
\Span\{H,\mathrm{ad}(\mathbb{D})H,\mathrm{ad}(\mathbb{D})^2H,\ldots,\mathrm{ad}(\mathbb{D})^{M -1}H\}
\]
contains the set $\{E^{i,j}\,:\, i,j\in\Z^{2}_{0,N},\, i\neq j\}$, which is easily seen to generate $\mathfrak{sl}_{\Z^2_{0,N}}(\R)$.

However in our setting the diagonal matrices $\mathbb{D}^k$ have an inversion symmetry $\mathbb{D}_{-i}^k = -\mathbb{D}^k_i$ and therefore there {\em cannot} be a matrix in $\mathfrak{h}$ with all differences distinct. Moreover, we do not have a matrix with {\em all} off diagonal entries non-zero due to the degeneracies present in $c_{j,k}$ and the presence of the Galerkin cut-off. Nevertheless, in \cite{BPS21} we are able to deduce the following sufficient condition on the family $\{\mathbb{D}^k\}$ that ensures \eqref{eq:Lie-suff-H} holds:
\begin{proposition}[Corollary 4.9 and Lemma 5.2 \cite{BPS21}]\label{prop:distinct}
Let $N\geq 8$. If for each $(i,j,\ell,m)\in (\Z^2_{0,N})^4$ satisfying $i+j+\ell+m = 0$ and
$(i+j,\ell+m)\neq 0,\, (i+\ell,j+m)\neq 0,\, (i+m,j+\ell) \neq 0$, there exists a $k\in\Z^2_{0,N}$ such that
\begin{equation}\label{eq:Distinctness}
	\mathbb{D}^k_i + \mathbb{D}_j^k+ \mathbb{D}^k_\ell+ \mathbb{D}^k_m \neq 0,
\end{equation}
then \eqref{eq:Lie-suff-H} holds.
\end{proposition}
The proof of Proposition \ref{prop:distinct} is not straightforward. However, its proof uses some similar ideas as the proof of \eqref{eq:Distinctness} but is otherwise significantly easier, so we only discuss the latter.

\subsubsection{Verifying the distinctness condition using computational algebraic geometry}

The distinctness condition \eqref{eq:Distinctness} is not a simple one to verify. Indeed, ignoring the Galerkin cut-off $N$ for now, $\mathbb{D}^k_i$ are rational algebraic expressions in the variables $(i,k,r)$ (being comprised of products and sums of the coefficients $c_{j,k}$), and therefore proving \eqref{eq:Distinctness} amounts to showing that the family of {\em Diophantine equations}\footnote{At least considering $r=1$ or another fixed, rational number.} 
\begin{equation}\label{eq:Diophantine}
	\mathbb{D}^k_i + \mathbb{D}_j^k+ \mathbb{D}^k_\ell+ \mathbb{D}^k_m = 0, \quad \text{for each} \quad k\in \Z^2_{0,N}
\end{equation}
have {\em no solutions} $(i,j,\ell,m,r)$ satisfying the constraints of Proposition \ref{prop:distinct}. Due to the complexity of the expression for $\mathbb{D}^k_i$, there is little hope to verify such a result by hand (the resulting polynomials are degree $16$ in $9$ variables). However, if one extends each of the $9$ variables $(i,j,\ell,m,r) = (i_1,i_2,j_1,j_2,\ell_1,\ell_2,m_1,m_2,r)$ to the algebraically closed field $\C$, then \eqref{eq:Diophantine} along with $i+j+\ell+m= 0$ defines a polynomial ideal $I$ with an associated algebraic variety $\mathbb{V}(I)$ in $\C^9$. Such a high dimensional variety is rather complicated due to the inherent symmetries of the rational equation in \eqref{eq:Diophantine}, however its analysis is nonetheless amenable to techniques from algebraic geometry, particularly the strong Nullstellensatz and computer algorithms for computing Gr\"obner bases (see \cite{Cox} for a review of the algebraic geometry concepts).
Indeed, without the Galerkin cut-off (the formal infinite dimensional limit), in \cite{BPS21} we proved, by computing Gr\"obner bases in rational arithmetic using the F4 algorithm \cite{F1999} implemented in the computer algebra system Maple \cite{maple}, that the identity $\mathbb{V}(I) = \mathbb{V}(g)$ holds, where $g$ is the following ``saturating'' polynomial
\[
	g(i,j,\ell,m,r) = r^2|i|_r^2|j|_r^2|\ell|_r^2|m|_r^2(|i+j|_r^2+|\ell+m|_r^2)(|i+\ell|_r^2 + |j+m|_r^2)(|i+m|^2_r + |j+\ell|_r^2)
\]
whose non-vanishing encodes the constraints in Proposition \ref{prop:distinct}, thereby showing that \eqref{eq:Distinctness} holds. 

Dealing with the Galerkin truncation adds significant difficulties to the proof as the associated rational system \eqref{eq:Diophantine} is instead piecewise defined (depending on $k$ and $N$) and therefore doesn't easily reduce to a problem about polynomial inconsistency. Nonetheless, by considering 34 different polynomial ideals associated to different possible algebraic forms, in \cite{BPS21} we were able to show that if $N$ is taken large enough (bigger than $392$ to be precise) then \eqref{eq:Distinctness} still holds with the Galerkin truncation present and therefore Theorem \ref{thm:2d-GSNS-Lie-gen} holds.

Finally, it is worth remarking that even without the Galerkin cut-off, the system of rational equations \eqref{eq:Diophantine} is complex enough to become computationally intractable (even for modern computer algebra algorithms) without some carefully chosen simplifications, variable orderings, choice of saturating polynomial $g$ and sheer luck; see \cite{BPS21} for more details.

%% file: LagChaos.tex
At present, the results above based on Proposition \ref{prop:fisherInfo} are restricted to finite dimensional problems. 
Indeed, even while the Fisher information can potentially be extended to infinite dimensions under certain conditions\footnote{If $X^\ast \nu \ll \nu$ and we define $\beta_{X^\ast}^\nu := \frac{\dee \tilde{X}^\ast \nu}{\dee \nu}$, then $FI(f) = \frac{1}{2} \sum_{k} \norm{\beta_{\tilde{X}^\ast}^\nu}_{L^2(\nu)}^2$, and there is no explicit dependence on any reference measure or Riemannian metric; see \cite{BBPS20} for more details.}, for any parabolic SPDE problem, we will always have $\lambda_\Sigma = -\infty$.
The existence of positive Lyapunov exponents for the infinite-dimensional, stochastic Navier-Stokes equations remains open as of the writing of this note.

However, there is another important problem in fluid mechanics where we have been able to make progress.
Consider the (infinite-dimensional) 2d Navier-Stokes equations\footnote{The 3D Navier-Stokes equations can be treated provided the $-\nu \Delta u_t$ is replaced with the hyperviscous damping $\nu \Delta^2 u_t$.} in $\mathbb T^2$, 
\begin{equation}\label{eq:2dSNS}
  \partial_t u_t  + \left( u_t \cdot \grad u_t + \grad p  - \nu \Delta u_t\right)  = \sum_k q_k e_k \dot{W}_t^{k},
  \quad \Div u_t  = 0, 
\end{equation}
where the $q_k \in \R$ and $e_k$ are eigenfunctions of the Stokes operator.
The \emph{Lagrangian flow map} $\varphi^t_{\omega,u} : \T^2 \mapsto \T^2$ is defined by the trajectories of particles moving with the fluid  
\begin{equation}
  \frac{\dee}{\dt}\varphi^t_{\omega,u}(x)  = u_t(\varphi^t_{\omega,u}(x)), \quad \varphi^0_{\omega,u}(x) = x, 
\end{equation}
where note that the diffeomorphism $\varphi^t_{\omega,u}$ depends on the initial velocity $u$ and the noise path $\omega$ and is therefore a co-cycle over the skew product $\Theta_t : \Omega\times H^s\circlearrowleft$, where $\Theta_t(\omega,u) = (\theta_t\omega,\Psi^t_{\omega}(u))$  and $\Psi^t_\omega: H^s \circlearrowleft$ is the 2d Navier Stokes flow on $H^s$ associated with \eqref{eq:2dSNS}. 
One can naturally ask whether or not $(u_t)$ is chaotic, as we have done in previous sections, or if the motion of particles immersed in the fluid is chaotic, e.g., if the Lagrangian Lyapunov exponent is strictly positive. 
The latter is known as \emph{Lagrangian chaos} \cite{BMOV05,CrisantiEtAl1991,AmonEtAl96,GalluccioVulpiani94,BalkovskyFouxon99,AntonsenEtAl96,YuanEtAl00} (to distinguish it from chaos of $(u_t)$ itself, which is sometimes called \emph{Eulerian chaos}).
While both are expected to be observed in turbulent flows, Lagrangian chaos is not incompatible
with Eulerian ``order'', i.e., a negative exponent for the $(u_t)$ process. 

In \cite{BBPS18} we proved, under the condition that $\abs{q_k} \approx \abs{k}^{-\alpha}$ for some $\alpha > 10$, that $\exists \lambda_1 > 0$ deterministic and independent of initial $x$ and initial velocity $u$ such that the following limit holds almost-surely: 
\begin{align}
\lim_{t \to \infty}\frac{1}{t} \log \abs{D_x \varphi^t_{\omega,u}} = \lambda_1 > 0. \label{def:LC} 
\end{align}
This Lagrangian chaos was later upgraded in \cite{BBPS19I,BBPS19II} to the much stronger property of uniform-in-diffusivity, almost-sure exponential mixing.
To formulate this notion, 
we consider $(g_t)$ a passive scalar solving the (random) advection-diffusion equation
\begin{equation}
  \partial_t g_t + u_t \cdot \grad g_t = \kappa \Delta g_t,\quad g_0 = g, 
\end{equation}
for $\kappa \in [0,1]$ and a fixed, mean-zero scalar $g \in L^2(\T^2)$. 
In \cite{BBPS19I,BBPS19II}, we proved that there exists a (deterministic) constant $\mu > 0$  such that for all $\kappa \in [0,1]$ and initial divergence free $u \in H^s$ (for some sufficiently large $s$), there exists a random constant $D = D(\omega,\kappa,u)$ such that for all $g \in H^1$ (mean-zero) 
\begin{equation}
\norm{g_t}_{H^{-1}} \leq D e^{-\mu t} \norm{g}_{H^1}
\end{equation}
where $D$ is almost-surely finite and satisfies the \emph{uniform-in-$\kappa$} moment bound (for some fixed constant $q$ and for any $\eta > 0$), 
\begin{equation}
\EE D^2 \lesssim_\eta (1 + \norm{u}_{H^s})^q e^{\eta \norm{u}_{H^1}^2}. 
\end{equation}
One can show that this result  is  essentially optimal up to getting sharper quantitative estimates on $\mu$ and $D$, at least if $\kappa = 0$ \cite{BBPS19I,BBPS19II}. This uniform, exponential mixing plays the key role in obtaining a proof of Batchelor's power spectrum \cite{Batchelor59} of passive scalar turbulence in some regimes \cite{BBPSbatch19}.  

Let us simply comment on the Lagrangian chaos statement \eqref{def:LC}, as it is most closely related to the rest of this note.
The main step is to deduce an analogue of Theorem \ref{cor:furstCriterion} for the Lagrangian flow map, using that while the Lagrangian flow map depends on an infinite dimensional Markov process, the Jacobian $D_x \varphi^t_{\omega,u}$ itself is finite dimensional.
This is done in our work \cite{BBPS18} by extending Furstenberg's criterion 
to handle general linear cocycles over infinite dimensional processes in the same way that
$D_x \varphi^t_{\omega, u}$ depends on the sample paths $(u_t)$. 


The Lagrangian flow is divergence-free, and thus the Lagrangian Lyapunov exponents satisfy $\lambda_\Sigma =0$ and $\lambda_1 \geq 0$, so ruling out the degenerate situations in 
Theorem \ref{cor:furstCriterion} would immediately imply $\lambda_1 > 0$. 
A key difficulty in this infinite-dimensional context is to ensure that the rigid invariant structures (now functions of the fluid velocity field $u$ and the Lagrangian tracer position $x$) in our analogue of Theorem \ref{cor:furstCriterion} vary \emph{continuously} as functions of $u$ and $x$. 
It is at this step that we require the non-degeneracy type condition on the noise $\abs{q_k} \gtrsim \abs{k}^{-\alpha}$, which is used to ensure that the Markov process $(u_t, \varphi^t(x))$ is strong Feller. 

At the time of writing, it remains an interesting open problem to extend our works \cite{BBPS18,BBPS19I,BBPS19II} to degenerate noise such as that used in \cite{HM06} or \cite{KNS20}. It bears remarking that the methods of \cite{HM06} apply to the one-point process $(u_t,\varphi^t(x))$ (this is used in our work \cite{BBPS19I}), however, it is nevertheless unclear how to prove Lagrangian chaos without a sufficiently strong analogue of Theorem \ref{cor:furstCriterion}, and it is unclear how to obtain such a theorem without the use of the strong Feller property.

%% file: forward.tex
The work we reviewed here raises a number of potential research directions. 



\medskip
\noindent {\bf Tighter hypoelliptic regularity estimates. }
The  scaling $\lambda_1^\eps \gg \eps$ that naturally follows from our above analysis is surely suboptimal -- even if the deterministic problem were to be completely integrable, the scaling is likely to be $O(\eps^{\gamma})$ for some $\gamma < 1$ depending on dimension (see e.g. \cite{pinsky1988lyapunov,baxendale2002lyapunov}).   
To begin with, one may attempt to strengthen the hypoelliptic regularity estimate by refining 
	the $\epsilon$ scaling to something like
	\[
	\| f^\eps \|_{W^{s, 1}}^2 \lesssim 1 + \frac{n \lambda_1^\eps - 2 \lambda_\Sigma^\eps}{\eps^{\gamma}},
	\]
for some constant $0 < \gamma < 1$. If such an estimate were true, the same compactness-rigidity argument 
of Theorem \ref{thm:critForEulerLikeIntro} would imply a scaling like $\lambda^\eps_1 \gtrsim \eps^{\gamma}$.
An improvement of this type seems plausible given the proof of Theorem \ref{thm:Hyp}.
It might be necessary, in general, to use a more specialized norm on the left-hand side, but local weak $L^1$ compactness, i.e. equi-integrability, is all that is really required for the compactness-rigidity argument to apply. 

\medskip
\noindent {\bf Beyond compactness-rigidity.}
Compactness-rigidity arguments may remain limited in their ability to obtain optimal or nearly optimal scalings for $\lambda_1$, regardless of the ways one can improve Theorem \ref{thm:Hyp}.
Another approach is to find some way to work more directly on $\eps > 0$. This was essentially
the approach of works \cite{pinsky1988lyapunov, baxendale2002lyapunov}, however, the method of these papers only applies if one has a nearly-complete understanding of the pathwise random dynamics.
We are unlikely to ever obtain such a complete understanding of the dynamics of models such as L96 or Galerkin-Navier-Stokes, but there may be hope that partial information, such as the isolation of robust, finite-time exponential growth mechanisms, could be used to obtain better lower bounds on $\norm {f^\eps}_{W^{s_*, 1}}$.
An approach with a vaguely related flavor for random perturbations of discrete-time systems, including the Chirikov standard map, was
carried out in the previous work \cite{blumenthal2017lyapunov}. 

\medskip
\noindent {\bf Finer dynamical information: moment Lyapunov exponents. }
Lyapunov exponents provide asymptotic exponential growth rates of the Jacobian, but they provide no quantitative information on how long it takes for this growth to be realized with high probability.
One tool to analyze this is the study of large deviations of the convergence of the sequences $\frac1t \log | D_x \Phi^t_\omega v|$.
The associated rate function  is the Legendre transform of the \emph{moment Lyapunov exponent function} $p \mapsto \Lambda(p) := \lim_{t \to \infty} \frac1t \log \E |D_x \Phi^t_\omega v|^p$ (the limit defining $\Lambda(p)$ exists and is independent of $(x, v, \omega)$ under fairly general conditions \cite{arnold1984formula}).
It would be highly interesting to see if the quantitative estimates obtained by e.g. Theorem \ref{thm:critForEulerLikeIntro} extend also to quantitative estimates on the moment Lyapunov exponents. 
We remark that the moment Lagrangian Lyapunov exponents play a key role in our works \cite{BBPS19I,BBPS19II}.

\medskip
\noindent {\bf Lyapunov times of small-noise perturbations of completely integrable systems}
The phase space of a completely integrable Hamiltonian flow is foliated by invariant torii along
which the dynamics is a translation flow-- such systems are highly ordered and non-chaotic. 
On the other hand, small perturbations of the Hamiltonian are 
known to break the most ``resonant'' of these torii, 
while torii with sufficiently ``non-resonant'' frequencies persist due to KAM theory. 
It is an interesting and highly challenging open problem to prove that this `breakage' 
results in the formation of a positive-volume set admitting a positive Lyapunov exponent. 
For the most part such problems are wide open, and related to the standard map conjecture
discussed in Section \ref{subsec:LEchallenges}.  
The recent work of Berger and Turaev \cite{berger2019herman} established a
 renormalization technique
for proving the \emph{existence} of smooth perturbations resulting in a positive Lyapunov exponent, but it remains open to affirm how `generic' such perturbations actually are. 

The following is a closely related \emph{stochastic} dynamics problem: starting from a 
completely integrable system and adding a small amount of noise, how many Lyapunov times
elapse for the random dynamics before the ``stochastic divergence'' timescale 
when the deterministic flow and the stochastic flow differ by $O(1)$?
Estimating the stochastic divergence timescale is essentially a large deviations problem, and 
has already been carried out for small random perturbations of completely integrable systems; see, e.g., \cite{freuidlin1994random}. On the other hand, estimating Lyapunov times beyond the crude $(\lambda_1^\eps)^{-1}$ estimate is a large deviations estimate for the convergence of finite-time Lyapunov exponents to their asymptotic value $\lambda_1^\eps$. The associated rate function in this case is the Legendre transform of the moment Lyapunov exponent $\Lambda(p)$ mentioned earlier; a positive result for the program described above would require quantitative-in-$\epsilon$ estimates on $\Lambda(p)$. 

\medskip
\noindent {\bf More general noise models. }
One simple potential extension is Theorem \ref{thm:critForEulerLikeIntro} to different types of multiplicative noise.
Another important extension would be to noise models which are not white-in-time, for example noise of the type used in \cite{KNS20}, which is challenging because our work is deeply tied to the elliptic nature of the generator $\tilde{\mathcal{L}}^\ast$. 
A simpler example of non-white forcing can be constructed from `towers' of coupled Ornstein-Uhlenbeck processes, which can be built to be $C^k$ in time for any $k \geq 0$ (see, e.g., \cite{BBPS18,BBPS19I} for details).

\medskip
\noindent {\bf Lagrangian chaos.}
 There are several directions of research to extend our results in \cite{BBPS18,BBPS19I,BBPS19II}, such as studying degenerate noise as in \cite{HM06,KNS20}, extending to more realistic physical settings such as bounded domains with stochastic boundary driving, and extending Proposition \ref{prop:fisherInfo} to the Lagrangian flow map in a variety of settings, which would help to facilitate quantitative estimates (note one will have to use the conditional density version so that one does not see the effect of the $\lambda_\Sigma$ associated to the Navier-Stokes equations themselves).